# EXTREMES ON TREES

By Tailen Hsing[1] and Holger Rootzén[2]

*Ohio State University and Chalmers University of Technology*

This paper considers the asymptotic distribution of the longest edge of the minimal spanning tree and nearest neighbor graph on $\mathbf{X}_1, \ldots, \mathbf{X}_{N_n}$ where $\mathbf{X}_1, \mathbf{X}_2, \ldots$ are i.i.d. in $\Re^2$ with distribution $F$ and $N_n$ is independent of the $\mathbf{X}_i$ and satisfies $N_n/n \to_p 1$. A new approach based on spatial blocking and a locally orthogonal coordinate system is developed to treat cases for which $F$ has unbounded support. The general results are applied to a number of special cases, including elliptically contoured distributions, distributions with independent Weibull-like margins and distributions with parallel level curves.

**1. Introduction.** Recall that the (Euclidean) minimal spanning tree (MST) on a finite set of points $(\mathbf{X}_1, \mathbf{X}_2, \ldots, \mathbf{X}_N)$ in $\Re^2$ is the connected graph with these points as vertices and with the minimum total edge length. The (Euclidean) nearest neighbor graph (NNG) on $(\mathbf{X}_1, \mathbf{X}_2, \ldots, \mathbf{X}_N)$ is the graph on which each point $\mathbf{X}_i$ is connected to its nearest neighbor in the set. In this paper the $\mathbf{X}_i$ are assumed to be random and we are interested in the asymptotic distribution of the longest edge on these graphs as $N \to \infty$.

Penrose (1997, 1998) considered these problems by assuming that the $\mathbf{X}_i$ are uniformly distributed in a unit cube or symmetrically normally distributed in $\Re^d$. The essential ideas are that (a) the lengths of the edges at any location in space depend primarily on the $\mathbf{X}_i$ in the vicinity of that location and as a result are asymptotically independent of the edges in other parts of the space and (b) the presence of an extremely long edge is a rare event and hence the likelihood of having multiple extremely long edges at any location is asymptotically negligible compared with the likelihood of having one such

Received June 2002; revised March 2004.

[1]Supported by the Texas Advanced Research Program, the Alexander von Humboldt Foundation and the German Science Foundation (Deutsche Forschungsgemeinschaft), SF 386, Statistical Analysis of Discrete Structures.

[2]Supported in part by the Knut and Alice Wallenberg Foundation.

edge there. Clearly (a) is also essential in proving central limit theorems for the total edge lengths; see Kesten and Lee (1996), Lee (1997) and Penrose (2000). In view of (a) and (b), the asymptotic distribution of the longest edge of MST or NNG can be established through a Poisson convergence of the number of extreme edges, which in Penrose (1997, 1998) is achieved through the Chen–Stein method.

Specifically, we will consider the case where the random graphs are based on $\mathbf{X}_1, \ldots, \mathbf{X}_{N_n}$ where $\mathbf{X}_1, \mathbf{X}_2, \ldots$ are i.i.d. with distribution $F$ and $N_n$ is independent of the $\mathbf{X}_i$ with $N_n/n \xrightarrow{p} 1$. We are primarily interested in the case where $F$ has an unbounded support (although our methodology works rather generally also if $F$ has bounded support). In particular, we will focus on those $F$ whose density is of the form

$$f(\mathbf{x}) = e^{-U(\mathbf{x})},$$

where $U(\mathbf{x})$ is regularly varying in some sense and satisfies suitable regularity conditions. This covers many elliptically contoured distributions, and in particular correlated normal distributions and distributions with independent Weibull-like marginals as special cases but also large classes of other distributions. Poisson approximation is the key idea. However, we use a direct approach of spatial blocking as opposed to the Chen–Stein method. Computations of integrals of the type $\int_A e^{-nF(S(\mathbf{x};r))} dF(\mathbf{x})$, where $S(\mathbf{x};r) = \{\mathbf{y} : |\mathbf{y} - \mathbf{x}| \leq r\}$ is the sphere centered at $\mathbf{x}$ with radius $r$, are a key part of the solution for the problem on hand. One of the novelties of our approach is the introduction of a locally orthogonal coordinate system with respect to the level curves of $U$, which enables particularly effective handling of such integrals.

This paper is structured as follows. Section 2 introduces the notation and a spatial blocking argument, as well as other preliminaries. The development of a locally orthogonal system is made in Section 3. The main results are given in Section 4. Sections 5 and 6 consider homogeneous level curves and parallel level curves, respectively, and most of the proofs are given in Sections 7 and 8.

Possible extensions of the results in this paper include (a) allowing the dimension $d$ to be general, (b) allowing the distance measure to be more general (e.g., considering weighted edges on the graphs), (c) considering the $k$-nearest neighbor graph in which each point is connected to its $k$ nearest neighbors. The solutions for these involve additional technical details but probably few new important ideas.

**2. Fundamentals.** For convenience of notation, denote by $\mathrm{MST}(N)$ and $\mathrm{NNG}(N)$, respectively, the MST and NNG on two-dimensional random variables $\mathbf{X}_1, \ldots, \mathbf{X}_N$ for any random variable $N$ defined on the same space as the $\mathbf{X}_i$. Also let $M_G$ be the longest edge of $G$, where in this paper $G$ will be either a MST or NNG.

As outlined in Section 1, the graphs of interest will be based on the points $\mathbf{X}_1, \ldots, \mathbf{X}_{N_n}$ where $\mathbf{X}_1, \mathbf{X}_2, \ldots$ are i.i.d. with distribution $F$, $N_n$ is positive, integer-valued with

$$N_n/n \xrightarrow{p} 1,$$

and $N_n$ and the $\mathbf{X}_i$ are independent. This assumption will hold throughout the paper. The random quantities whose asymptotic distributions we study in this paper are $M_{\mathrm{MST}(N_n)}$ and $M_{\mathrm{NNG}(N_n)}$. We will first assume that $N_n$ is Poisson distributed with mean $n$, in which case the points $\mathbf{X}_1, \ldots, \mathbf{X}_{N_n}$ can be thought of as the points of a Poisson process with intensity measure $nF$. The independent-increment property of the Poisson process offers an obvious advantage in proving Poisson convergence. We will later show how the result based on the Poisson assumption can be extended to the general class of point processes described here.

Temporal blocking is a common technique for proving limit theorems for weakly dependent random variables. See Ibragimov and Linnik (1971) and Leadbetter, Lindgren and Rootzén (1983). Our first theorem, which basically is a Poisson convergence result for $M_{\mathrm{NNG}(N_n)}$, for the case $N_n \sim \mathrm{Poisson}(n)$, is based on a spatial blocking argument. For $r > 0$ and measurable set $A$, define

$$(1) \qquad \mu_n^{(1)}(A, r) = n \int_A e^{-nF(S(\mathbf{x};r))} \, dF(\mathbf{x}).$$

Now, let us say that points with their nearest neighbor at least $r$ away are *r-separate*. With this terminology, $\mu_n^{(1)}(A, r)$ is the expected number of $r$-separate points in $A$. Similarly, let

$$\mu_n^{(2)}(A, r) = n^2 \int_A \int_A I_{(r < |\mathbf{x}-\mathbf{y}| \leq 2r)} e^{-nF(S(\mathbf{x};r) \cup S(\mathbf{y};r))} \, dF(\mathbf{x}) \, dF(\mathbf{y})$$

be the expected number of pairs of $r$-separate points such that the distance between the points is larger than $r$ and smaller than $2r$. Here and elsewhere, $|\cdot|$ denotes the Euclidean norm.

In the theorem, for each $n$, $A_{n,i}, 1 \leq i \leq k_n$, are suitably large and suitably separated spatial blocks. The separation has to be large enough to make the occurrences of $r_n$-separate points independent from block to block. This is ensured by condition (a). Condition (c) says that nothing of importance happens on the leftover parts between the blocks. Uniform asymptotic negligibility of the number of $r_n$-separate points in the individual blocks follows from (d). Condition (e) prevents clustering of $r_n$-separate points. Finally, condition (b) is the basic norming condition of convergence of the expected number of $r_n$-separate points to the mean of the limiting Poisson distribution.

THEOREM 1. *Let $N_n$ have a Poisson distribution with mean $n$ and let $\{k_n\}$ be a sequence of positive constants tending to $\infty$. Suppose that $\{r_n\}$ is a sequence of positive constants and for each $n$, $A_{n,1}, \ldots, A_{n,k_n}$ are disjoint measurable sets in $\Re^2$ such that:*

(a) $\min_{1 \leq i \neq j \leq k_n} \inf(|\mathbf{x} - \mathbf{y}| : \mathbf{x} \in A_{n,i}, \mathbf{y} \in A_{n,j}) > 2r_n$ *for all large $n$,*
(b) $\lim_{n \to \infty} \mu_n^{(1)}(\Re^2, r_n) = $ *some* $\tau \in (0, \infty)$,
(c) $\lim_{n \to \infty} \mu_n^{(1)}((\bigcup_{i=1}^{k_n} A_{n,i})^c, r_n) = 0$,
(d) $\lim_{n \to \infty} \max_{1 \leq i \leq k_n} \mu_n^{(1)}(A_{n,i}, r_n) = 0$,
(e) $\lim_{n \to \infty} \max_{1 \leq i \leq k_n} [\mu_n^{(2)}(A_{n,i}, r_n) / \mu_n^{(1)}(A_{n,i}, r_n)] = 0$.

*Then the number of $r_n$-separate points asymptotically has a Poisson distribution with mean $\tau$, and thus in particular,*

$$P(M_{\mathrm{NNG}(N_n)} \leq r_n) \to e^{-\tau}. \tag{2}$$

PROOF. For convenience write $\mathcal{P}_n = \{\mathbf{X}_1, \ldots, \mathbf{X}_{N_n}\}$. For any set $A$, define

$$N(\mathbf{x}, A) = \inf(|\mathbf{x} - \mathbf{y}| : \mathbf{y} \in A \setminus \{\mathbf{x}\}),$$

namely, the distance from $\mathbf{x}$ to its nearest neighbor in $A$. We will show the stronger result that $\sum_{\mathbf{x} \in \mathcal{P}_n} I_{(N(\mathbf{x}, \mathcal{P}_n) > r_n)}$ converges in distribution to the Poisson distribution with mean $\tau$, from which (2) follows at once. Let

$$\tilde{\mu}_n^{(2)}(A, r) = n^2 \int_A \int_A I_{(|\mathbf{x} - \mathbf{y}| > r)} e^{-nF(S(\mathbf{x};r) \cup S(\mathbf{y};r))} \, dF(\mathbf{x}) \, dF(\mathbf{y})$$

be the expected number of pairs of $r$-separate points in $A$ with the points at least a distance $r$ apart. It is easy to check that

$$E\left(\sum_{\mathbf{x} \in \mathcal{P}_n \cap A} I_{(N(\mathbf{x}, \mathcal{P}_n) > r)}\right) = \mu_n^{(1)}(A, r) \tag{3}$$

and

$$E\left(\sum_{\mathbf{x}, \mathbf{y} \in \mathcal{P}_n \cap A,\, \mathbf{x} \neq \mathbf{y}} I_{(N(\mathbf{x}, \mathcal{P}_n) \wedge N(\mathbf{y}, \mathcal{P}_n) > r)}\right) = \tilde{\mu}_n^{(2)}(A, r). \tag{4}$$

Write

$$\sum_{\mathbf{x} \in \mathcal{P}_n} I_{(N(\mathbf{x}, \mathcal{P}_n) > r_n)} = \sum_{i=1}^{k_n} \sum_{\mathbf{x} \in \mathcal{P}_n \cap A_{n,i}} I_{(N(\mathbf{x}, \mathcal{P}_n) > r_n)} + \sum_{\mathbf{x} \in \mathcal{P}_n \cap (\bigcup_{i=1}^{k_n} A_{n,i})^c} I_{(N(\mathbf{x}, \mathcal{P}_n) > r_n)}.$$

Roughly speaking, in the following we will establish that the second sum is negligible and for the first sum that its expectation tends to $\tau$, the summands are infinitesimal and the probability that any summand is bigger than 1 is negligible compared with that of it being bigger than 0. More precisely, we will show that

$$\sum_{\mathbf{x}\in\mathcal{P}_n\cap A_{n,i}} I_{(N(\mathbf{x},\mathcal{P}_n)>r_n)}, \qquad 1\leq i\leq k_n,$$

(5)

are independent random variables for each $n$,

(6) $$P\left(\sum_{\mathbf{x}\in\mathcal{P}_n\cap(\bigcup_{i=1}^{k_n} A_{n,i})^c} I_{(N(\mathbf{x},\mathcal{P}_n)>r_n)} > 0\right) \to 0,$$

(7) $$\max_{1\leq i\leq k_n} P\left(\sum_{\mathbf{x}\in\mathcal{P}_n\cap A_{n,i}} I_{(N(\mathbf{x},\mathcal{P}_n)>r_n)} > 0\right) \to 0,$$

(8) $$\sum_{i=1}^{k_n} P\left(\sum_{\mathbf{x}\in\mathcal{P}_n\cap A_{n,i}} I_{(N(\mathbf{x},\mathcal{P}_n)>r_n)} > 0\right) \to \tau,$$

(9) $$\max_{1\leq i\leq k_n} \frac{P(\sum_{\mathbf{x}\in\mathcal{P}_n\cap A_{n,i}} I_{(N(\mathbf{x},\mathcal{P}_n)>r_n)} > 1)}{P(\sum_{\mathbf{x}\in\mathcal{P}_n\cap A_{n,i}} I_{(N(\mathbf{x},\mathcal{P}_n)>r_n)} > 0)} \to 0.$$

It is straightforward to see [cf. Corollary 7.5 of Kallenberg (1983)] that (5)–(9) imply that $\sum_{\mathbf{x}\in\mathcal{P}_n} I_{(N(\mathbf{x},\mathcal{P}_n)>r_n)}$ converges in distribution to the Poisson distribution with mean $\tau$.

We now show (5)–(9). First, it is clear that $\sum_{\mathbf{x}\in\mathcal{P}_n\cap A_{n,i}} I_{(N(\mathbf{x},\mathcal{P}_n)>r_n)}$ is completely determined by the set $\mathcal{P}_n \cap (A_{n,i})_{r_n}$ where $(A)_\delta = \{\mathbf{x}: d(\mathbf{x},A)\leq \delta\}$. Hence, by the independent-increment property of the Poisson process and condition (a) of Theorem 1, we conclude that (5) holds. To show (6)–(9), we first note the Bonferroni inequality:

(10) $$\mu_n^{(1)}(A,r_n) - \tilde{\mu}_n^{(2)}(A,r_n) \leq P\left(\sum_{\mathbf{x}\in\mathcal{P}_n\cap A} I_{(N(\mathbf{x},\mathcal{P}_n)>r_n)} > 0\right) \leq \mu_n^{(1)}(A,r_n).$$

Using the rightmost inequality with $A = (\bigcup_{i=1}^{k_n} A_{n,i})^c$ and $A = A_{n,i}$, respectively, (6) and (7) follow from conditions (c) and (d), respectively.

Next, note that since $F(S(\mathbf{x};r_n) \cup S(\mathbf{y};r_n)) = F(S(\mathbf{x};r_n)) + F(S(\mathbf{y};r_n))$ for $|\mathbf{x}-\mathbf{y}|>2r_n$,

$$\tilde{\mu}_n^{(2)}(A,r_n) - \mu_n^{(2)}(A,r_n) = n^2 \int_A \int_A I_{(|\mathbf{x}-\mathbf{y}|>2r_n)} e^{-nF(S(\mathbf{x};r_n)\cup S(\mathbf{y};r_n))} \, dF(\mathbf{x}) \, dF(\mathbf{y})$$

$$\leq (\mu_n^{(1)}(A,r_n))^2.$$

Hence, by conditions (d) and (e),

(11) $$\lim_{n\to\infty} \max_{1\leq i\leq k_n} \frac{\tilde{\mu}_n^{(2)}(A_{n,i},r_n)}{\mu_n^{(1)}(A_{n,i},r_n)} = 0.$$

By (10) and (11), conditions (b) and (c) now imply

$$\lim_{n\to\infty} \sum_{i=1}^{k_n} P\left(\sum_{\mathbf{x}\in\mathcal{P}_n\cap A_{n,i}} I_{(N(\mathbf{x},\mathcal{P}_n)>r_n)} > 0\right) = \lim_{n\to\infty} \sum_{i=1}^{k_n} \mu_n^{(1)}(A_{n,i}, r_n)$$

$$= \lim_{n\to\infty} \mu_n^{(1)}(\Re, r_n) = \tau.$$

Finally, since

$$P\left(\sum_{\mathbf{x}\in\mathcal{P}_n\cap A} I_{(N(\mathbf{x},\mathcal{P}_n)>r_n)} > 1\right) \le P\left(\sum_{\mathbf{x},\mathbf{y}\in\mathcal{P}_n\cap A, \mathbf{x}\neq\mathbf{y}} I_{(N(\mathbf{x},\mathcal{P}_n)\wedge N(\mathbf{y},\mathcal{P}_n)>r_n)} > 0\right)$$

$$\le \tilde{\mu}_n^{(2)}(A, r_n),$$

(9) follows from (10) and (11). □

As a simple example, if $F$ is the uniform distribution on $[0,1]\times[0,1]$, then it is straightforward to verify (a)–(e) of Theorem 1 by taking $r_n = (\frac{\log n - \log \tau}{n\pi})^{1/2}$ [cf. (1) of Penrose (1997)] and the $A_{n,i}$ to be the sets $[(j-1)n^{-1/4}, jn^{-1/4} - n^{-1/3}] \times [(k-1)n^{-1/4}, kn^{-1/4} - n^{-1/3}]$, $j, k = 1, \ldots, [n^{1/4}]$. The advantage of this approach will be more obvious for more complicated distributions.

The next result shows how to generalize Theorem 1 by removing the Poisson assumption on $N_n$.

THEOREM 2. *Suppose that*

$$\sup_n n \int \bar{F}^n(S(\mathbf{x}; r_n)) \, dF(\mathbf{x}) < \infty \tag{12}$$

*and for some $\delta < 1$,*

$$\sup_n n^2 \int F(S(\mathbf{x}; r_n)) \bar{F}^{\delta n}(S(\mathbf{x}; r_n)) \, dF(\mathbf{x}) < \infty, \tag{13}$$

*where $\bar{F} = 1 - F$. Then*

$$P(M_{\mathrm{NNG}(N_n)} \le r_n < M_{\mathrm{NNG}(n)}) + P(M_{\mathrm{NNG}(n)} \le r_n < M_{\mathrm{NNG}(N_n)}) \to 0$$

*for any sequence of positive, integer-valued random variables $N_n$ with $N_n/n \xrightarrow{p} 1$.*

PROOF. First for any positive constant $\varepsilon_n$,

$$P(M_{\mathrm{NNG}(N_n)} \le r_n < M_{\mathrm{NNG}(n)}) + P(M_{\mathrm{NNG}(n)} \le r_n < M_{\mathrm{NNG}(N_n)})$$

$$\le P(|N_n - n| > n\varepsilon_n) + \sum_{|j-n|\le n\varepsilon_n} P(N_n = j)$$

$$\times \Big(P(M_{\mathrm{NNG}(j)} \le r_n < M_{\mathrm{NNG}(n)})$$

$$+ P(M_{\mathrm{NNG}(n)} \le r_n < M_{\mathrm{NNG}(j)})\Big).$$

Pick $\varepsilon_n$ to tend to 0 slowly enough so that the first term tends to 0 and so we only have to deal with the second term. It suffices to show that

$$(14) \qquad \max_{j \in [n(1-\varepsilon_n), n(1+\varepsilon_n)]} P(M_{\mathrm{NNG}(j)} \leq r_n < M_{\mathrm{NNG}(n)}) \to 0$$

and

$$(15) \qquad \max_{j \in [n(1-\varepsilon_n), n(1+\varepsilon_n)]} P(M_{\mathrm{NNG}(n)} \leq r_n < M_{\mathrm{NNG}(j)}) \to 0.$$

To show (14), suppose first that $j \in [n(1-\varepsilon_n), n-1]$. Observe that if $M_{\mathrm{NNG}(j)} \leq r_n < M_{\mathrm{NNG}(n)}$, then there exists $i \in \{j+1,\ldots,n-1\}$ such that $\bigwedge_{1 \leq s \leq n, s \neq i} |X_i - X_s| > r_n$. Hence

$$\max_{j \in [n(1-\varepsilon_n), n-1]} P(M_{\mathrm{NNG}(j)} \leq r_n < M_{\mathrm{NNG}(n)}) \leq n\varepsilon_n \int \bar{F}^{n-1}(S(\mathbf{x}; r_n))\, dF(\mathbf{x}) \to 0$$

by (12). Next, let $j \in [n+1, n(1+\varepsilon_n)]$. Suppose that $M_{\mathrm{NNG}(j)} \leq r_n < M_{\mathrm{NNG}(n)}$ and that the longest edge on $\mathrm{NNG}(n)$ initiates from $\mathbf{X}_k$ so that $\bigwedge_{1 \leq s \leq n, s \neq k} |\mathbf{X}_k - \mathbf{X}_s| > r_n$. Since the longest edge becomes $\leq r_n$ by adding points $\mathbf{X}_{n+1},\ldots,\mathbf{X}_j$ to the graph, one of those additional points must be within a distance of $r_n$ from $\mathbf{X}_k$. Thus,

$$P(M_{\mathrm{NNG}(j)} \leq r_n < M_{\mathrm{NNG}(n)})$$
$$\leq P\bigg(\text{for some } k \in \{1,\ldots,n\} \text{ and } \ell \in \{n+1,\ldots,j\},$$
$$\bigwedge_{1 \leq s \leq n, s \neq k} |X_k - X_s| > r_n \text{ but } |X_k - X_\ell| \leq r_n\bigg)$$
$$\leq n^2 \varepsilon_n \int F(S(\mathbf{x}; r_n)) \bar{F}^{n-1}(S(\mathbf{x}; r_n))\, dF(\mathbf{x}) \to 0$$

uniformly for $j \in [n+1, n(1+\varepsilon_n)]$ by (13). This proves (14). The proof of (15) is similar, where the main difference is that for $j \in [n(1-\varepsilon_n), n-1]$,

$$P(M_{\mathrm{NNG}(n)} \leq r_n < M_{\mathrm{NNG}(j)}) \leq n^2 \varepsilon_n \int F(S(\mathbf{x}; r_n)) \bar{F}^{j-1}(S(\mathbf{x}; r_n))\, dF(\mathbf{x}),$$

which, again by (13), tends to 0 uniformly for $j \in [n(1-\varepsilon_n), n-1]$. $\square$

We note that (12) follows readily from (b) of Theorem 1 since $\bar{F}(S(\mathbf{x}; r_n))^n \leq e^{-nF(S(\mathbf{x};r_n))}$. Indeed, for the cases that we will consider, the conditions (12) and (13) are both naturally satisfied. Thus, once we have the asymptotic distribution of the extreme edge length for the Poisson NNG, we can extend that at once to a more general class of NNGs.

The following result gives an argument for deriving the asymptotic distribution of $M_{\mathrm{MST}(N_n)}$ from that of $M_{\mathrm{NNG}(N_n)}$, in light of the fact that $M_{\mathrm{MST}(N_n)} \geq M_{\mathrm{NNG}(N_n)}$ (see proof below).

THEOREM 3. *Suppose that for any sequence of positive constants $m_n$ with $m_n/n \to 1$, we have*

$$(16) \quad n \int [\bar{F}^{m_n}(S(\mathbf{x}; r_n) \cap S(0; |\mathbf{x}|)) - \bar{F}^{m_n}(S(\mathbf{x}; r_n))] I(|\mathbf{x}| > r_n/2) \, dF(\mathbf{x}) \to 0.$$

*Then*

$$P(M_{\mathrm{NNG}(N_n)} \leq r_n < M_{\mathrm{MST}(N_n)}) \to 0$$

*for any sequence of positive, integer-valued random variables $N_n$ with $N_n/n \xrightarrow{p} 1$.*

PROOF. The proof is basically the same as that of Lemma 4 of Penrose (1998). For completeness, the essential ideas are reproduced here. Again write

$$\begin{aligned}
&P(M_{\mathrm{NNG}(N_n)} \leq r_n < M_{\mathrm{MST}(N_n)}) \\
(17) \quad &\leq P(|N_n/n - 1| > \varepsilon_n) \\
&\quad + \sum_{|m/n-1| \leq \varepsilon_n} P(N_n = m) P(M_{\mathrm{NNG}(m)} \leq r_n < M_{\mathrm{MST}(m)}),
\end{aligned}$$

where $\varepsilon_n$ tends to zero slowly enough so that the first term tends to 0. Now make $\mathbf{X}_1, \ldots, \mathbf{X}_m$ a graph by including an edge between each pair of points at a distance at most $r$. Denote the resulting graph by $G_r$. Clearly, for a small enough $r$, $G_r$ comprises connected subgraphs, called $r$-clusters, which are disconnected with one another. Observe that

$M_{\mathrm{MST}(m)} = \inf\{r : G_r \text{ is connected}\},$

$M_{\mathrm{NNG}(m)} = \inf\{r : G_r \text{ does not contain an } r\text{-cluster which is a singleton}\}.$

Thus, $M_{\mathrm{MST}(m)} \geq M_{\mathrm{NNG}(m)}$. Suppose $M_{\mathrm{NNG}(m)} \leq r_n < M_{\mathrm{MST}(m)}$. Then it means that $G_{r_n}$ is disconnected and each $r_n$-cluster has at least two points. Take an $r_n$-cluster and let $\mathbf{x}$ be the vertex in the cluster which is closest to 0. Write $\mathcal{I}_m^{\mathbf{x}} = \{\mathbf{X}_1, \ldots, \mathbf{X}_m\} \setminus \{\mathbf{x}\}$. Clearly

$$\mathcal{I}_m^{\mathbf{x}} \cap S(\mathbf{x}; r_n) \neq \varnothing \quad \text{and} \quad \mathcal{I}_m^{\mathbf{x}} \cap S(\mathbf{x}; r_n) \cap S(0; |\mathbf{x}|) = \varnothing$$

since $S(\mathbf{x}; r_n)$ contains points and only points belonging to the same $r_n$-cluster while $S(0; |\mathbf{x}|) \setminus \{\mathbf{x}\}$ does not contain points from the same $r_n$-cluster. This means that if $M_{\mathrm{NNG}(m)} \leq r_n < M_{\mathrm{MST}(m)}$, then there are at least two points belonging to different $r_n$-clusters having the described property, where one of them must have modulus bigger than $r_n/2$. Thus,

$$\begin{aligned}
&P(M_{\mathrm{NNG}(m)} \leq r_n < M_{\mathrm{MST}(m)}) \\
&\leq P\bigg(\sum_{i=1}^{m} I(|\mathbf{X}_i| > r_n/2, \mathcal{I}_m^{\mathbf{X}_i} \cap S(\mathbf{X}_i; r_n) \neq \varnothing
\end{aligned}$$

$$\text{and } \mathcal{I}_m^{\mathbf{X}_i} \cap S(\mathbf{X}_i; r_n) \cap S(0; |\mathbf{X}_i|) = \varnothing) \geq 1\bigg)$$

$$\leq E\bigg(\sum_{i=1}^m I(|\mathbf{X}_i| > r_n/2, \mathcal{I}_m^{\mathbf{X}_i} \cap S(\mathbf{X}_i; r_n) \neq \varnothing$$

$$\text{and } \mathcal{I}_m^{\mathbf{X}_i} \cap S(\mathbf{X}_i; r_n) \cap S(0; |\mathbf{X}_i|) = \varnothing)\bigg)$$

$$= E\bigg(\sum_{i=1}^m I(|\mathbf{X}_i| > r_n/2, \mathcal{I}_m^{\mathbf{X}_i} \cap S(\mathbf{X}_i; r_n) \cap S(0; |\mathbf{X}_i|) = \varnothing)\bigg)$$

$$- E\bigg(\sum_{i=1}^m I(|\mathbf{X}_i| > r_n/2, \mathcal{I}_m^{\mathbf{X}_i} \cap S(\mathbf{X}_i; r_n) = \varnothing)\bigg)$$

$$= m\int [\bar{F}^{m-1}(S(\mathbf{x}; r_n) \cap S(0; |\mathbf{x}|)) - \bar{F}^{m-1}(S(\mathbf{x}; r_n))]I(|\mathbf{x}| > r_n/2)\,dF(\mathbf{x}),$$

which tends to 0 by (16). The result follows from this in view of (17). $\square$

While it may not be easy to make an intuitive connection between the condition (16) and its conclusion in Theorem 3, (16) holds quite naturally for the cases that we study in this paper.

It might also be worth noting that it is easy to find examples where the Poisson convergence which is used for Theorem 1 does not hold. In fact, this is the rule rather than the exception in the one-dimensional case, and for example can be seen to be the case for the one-dimensional density $e^{-|x|}/2$. A simple two-dimensional example is then obtained by letting the two-dimensional distribution be concentrated on the $x_1$-axis and have this density. If this distribution is mixed with, say, the standard bivariate normal distribution, a more genuine two-dimensional example where Poisson convergence does not occur is obtained. This example can also be simply modified to have a smooth density.

**3. Time-varying, locally orthogonal coordinate system.** In the remaining part of the paper, we will focus on densities of the form

$$f(\mathbf{x}) = e^{-U(\mathbf{x})},$$

where $U(\mathbf{x})$ is continuous and each level curve $(U = u) := \{\mathbf{x}: U(\mathbf{x}) = u\}$ is a closed and convex curve. This guarantees that $f$ is monotonically decreasing in some sense. We first introduce a "locally orthogonal" coordinate system which is useful for our computations, in particular for computations of the basic quantity $\mu_n^{(1)}(A, r)$ defined by (1). Let

$$\nabla U(\mathbf{x}) = (U^{(1,0)}(\mathbf{x}), U^{(0,1)}(\mathbf{x}))$$

be the gradient of $U(\mathbf{x})$ at $\mathbf{x}$. This gradient is throughout assumed to be continuous.

For points $\mathbf{x}$ for which $U(\mathbf{x})$ is large, we define the transformation $\mathbf{x} \to (\ell, u)$ where $u = U(\mathbf{x})$ and $\ell$ is determined in the following manner. We use a level curve $(U = w)$, with $w$ specified below, as a reference. If $\mathbf{x}$ belongs to the reference curve $(U = w)$, let $\ell(\mathbf{x})$ be the arc length from an arbitrary but fixed point $\mathbf{x}_o$ to $\mathbf{x}$, measured counterclockwise, say, on the curve $(U = w)$. In general, if $\mathbf{x}$ is on the level curve $(U = u)$, where $u \neq w$ we let $\ell(\mathbf{x}) = \ell(\mathbf{x}')$ where $\mathbf{x}'$ lies on $(U(\mathbf{x}) = w)$ and $\mathbf{x}$ and $\mathbf{x}'$ are connected by a curve which is orthogonal to each level curve between $(U = w)$ and $(U = u)$, in the sense that the normal of any intersecting level curve at the point of intersection is in the same direction (modulo $\pi$) as the tangent on the connecting curve at the point. The assumption that $\nabla(\mathbf{x})$ is continuous assures that the transformation $\mathbf{x} \to (\ell, u)$ is one-to-one. The reason for this choice of the second coordinate, $\ell$, is that it gives the Jacobian for the transformation to the new coordinate system a relatively simple form. The important thing to keep in mind is that this coordinate system depends on the reference curve. In what follows we refer to the index $n$ in $N_n$ as "time." In the computations that follow, we see that at time $n$ when $n$ is large, all the action takes place in a neighborhood of the level curve $(U = \log n)$ (see proof of Theorem 4). Hence at time $n$ it is natural to set $w = \log n$ and use $(U = \log n)$ as the reference curve in defining $\ell$. For this reason we shall throughout the rest of this paper adhere to this particular coordinate system. In doing so, here and elsewhere, the reference of $\ell$ to time will be suppressed (i.e., instead of $\ell_n$ we simply write $\ell$) for convenience of notation. Also, whenever there is no ambiguity the notation $\ell, \mathbf{x}$ will denote the functions $\ell(\mathbf{x}), \mathbf{x}(\ell, u)$ as well as their possible values.

We now derive the Jacobian for the coordinate change from $\mathbf{x}$ to $(u, \ell)$. Note that the unit normal vector at $\mathbf{x}$ on the corresponding level curve is

$$\frac{\nabla U(\mathbf{x})}{|\nabla U(\mathbf{x})|}.$$

By the way in which the transformation is defined,

$$(18) \qquad \frac{\partial \mathbf{x}}{\partial u} = \frac{1}{|\nabla U(\mathbf{x})|} \frac{\nabla U(\mathbf{x})}{|\nabla U(\mathbf{x})|} = \frac{\nabla U(\mathbf{x})}{|\nabla U(\mathbf{x})|^2}.$$

We now derive $(\frac{\partial x_1}{\partial \ell}, \frac{\partial x_2}{\partial \ell})$. Clearly the unit tangent vector at $\mathbf{x}$ on the corresponding level curve is

$$\frac{1}{|\nabla U(\mathbf{x})|} (U^{(0,1)}(\mathbf{x}), -U^{(1,0)}(\mathbf{x})).$$

If **x** is such that $U(\mathbf{x}) = \log n$, then $\ell$ corresponds to the actual arc length and so

$$(19) \qquad \frac{\partial \mathbf{x}}{\partial \ell} = \text{ unit tangent at } \mathbf{x} = \frac{1}{|\nabla U(\mathbf{x})|}(U^{(0,1)}(\mathbf{x}), -U^{(1,0)}(\mathbf{x})),$$

in which case the Jacobian is equal to

$$\left|\frac{\partial \mathbf{x}}{\partial(\ell, u)}\right| = \left|\frac{\partial x_1}{\partial u}\frac{\partial x_2}{\partial \ell} - \frac{\partial x_1}{\partial \ell}\frac{\partial x_2}{\partial u}\right| = \frac{1}{|\nabla U(\mathbf{x})|}.$$

In general, the arc length from $\mathbf{x}(\ell_1, u)$ to $\mathbf{x}(\ell_2, u)$ is computed as

$$(20) \qquad \int_{\ell_1}^{\ell_2} |\mathbf{x}^{(1,0)}(\ell, u)| \, d\ell.$$

Thus, if $\ell$ changes by $d\ell$, then the arc length changes by $|\mathbf{x}^{(1,0)}(\ell, u)| \, d\ell$ and therefore

$$\frac{\partial \mathbf{x}}{\partial \ell} = \frac{|\mathbf{x}^{(1,0)}(\ell, u)|}{|\nabla U(\mathbf{x})|}(U^{(0,1)}(\mathbf{x}), -U^{(1,0)}(\mathbf{x})).$$

Consequently the Jacobian is

$$(21) \qquad \left|\frac{\partial \mathbf{x}}{\partial(\ell, u)}\right| = \frac{|\mathbf{x}^{(1,0)}(\ell, u)|}{|\nabla U(\mathbf{x}(\ell, u))|}.$$

For convenience write

$$(22) \qquad \xi(\ell, u) = |\nabla U(\mathbf{x}(\ell, u))|$$

and

$$(23) \qquad \lambda(u) = \text{ length of the level curve } (U = u).$$

Thus, changing variables $\mathbf{x} \to (\ell, u)$ gives

$$(24) \qquad \begin{aligned} \mu_n^{(1)}(A, r) &= \int_A e^{-nF(S(\mathbf{x};r))} \, dF(\mathbf{x}) \\ &= \int_{\tilde{A}} e^{-nF(S(\mathbf{x}(\ell,u);r))} e^{-u} \frac{|\mathbf{x}^{(1,0)}(\ell, u)|}{\xi(\ell, u)} \, d\ell \, du, \end{aligned}$$

since $dF(\mathbf{x}) = e^{-U(\mathbf{x})} \, d\mathbf{x}$, where $\tilde{A}$ is the image of $A$ under the transformation from $\mathbf{x}$ to $(\ell, u)$.

**4. Main results.** In this section we formulate conditions directly in terms of $U(\mathbf{x})$, the negative logarithm of the density $f(\mathbf{x}) = e^{-U(\mathbf{x})}$, which lead to convergence of the longest edges of the nearest neighbor graph and the minimal spanning tree. As stated in the beginning of Section 3, we focus on functions $U$ for which the gradient is continuous and the level curves are convex. Additional conditions on smoothness and other aspects of $U$ will be stated below in Assumptions A1–A6. The conditions are rather technical. Some discussion of their meaning is given after the statement of Theorem 5.

The core of the problem is to find sequences $r_n$ for which the mean number of $r_n$-separate points converges, that is, which satisfy $\lim_{n\to\infty} \mu_n^{(1)}(\Re^2, r_n) = $ some $\tau \in (0, \infty)$. Denote by $\ell_o = \ell_{n,o}$ any point for which

$$(25) \qquad \xi_n := \xi(\ell_o, \log n) = \min_\ell \xi(\ell, \log n).$$

We assume that there exist finite positive constants $c_1, c_2$ and a sequence $\eta_n$ which satisfy Assumption A5, for $r_n$ defined as

$$(26) \qquad r_n := \frac{\eta_n - c_1 \log \eta_n - \log(c_2^{-1} \tau \sqrt{2\pi})}{\xi_n}.$$

Throughout assume that $r_n$ is given by this expression.

For convenience write $\log_2 = \log\log$ and $\log_3 = \log\log\log$. The following set of assumptions are needed:

ASSUMPTION A1. For each large $u$, given any $\varepsilon > 0$ there exists $\delta > 0$ such that $|\nabla U(\mathbf{x})/|\nabla U(\mathbf{x})| - \nabla U(\mathbf{x}')/|\nabla U(\mathbf{x}')|| < \varepsilon$ for all $\mathbf{x}, \mathbf{x}'$ on $U = u$ with their distance on the curve less than $\delta \lambda(u)$.

ASSUMPTION A2. $\limsup_{u\to\infty} \frac{\sup_{\mathbf{x}: U(\mathbf{x})=u} |\nabla U(\mathbf{x})|}{\inf_{\mathbf{x}: U(\mathbf{x})=u} |\nabla U(\mathbf{x})|} < \infty.$

ASSUMPTION A3. For all large $u$, both $\lambda(u)$ and $\inf_{\mathbf{x}: U(\mathbf{x})=u} |\nabla U(\mathbf{x})|$ belong to the range $[u^{-\rho}, u^\rho]$ for some $\rho \in (0, \infty)$.

ASSUMPTION A4. $\lim_{U(\mathbf{x})\to\infty} \frac{|U^{(i,j)}(\mathbf{x})|}{|\nabla U(\mathbf{x})|^2}(\log U(\mathbf{x}))^2 = 0$ for any $i, j \geq 0$ with $i + j = 2$.

ASSUMPTION A5. There exist positive constants $c_1, c_2$ and a sequence of constants $\eta_n \to \infty$ with $\eta_n = O(\log_2 n)$ such that [with $\xi(\ell, \log n) = |\nabla U(\mathbf{x}(\ell, \log n))|$, as before]:

(a) $\lim_{n\to\infty} e^{-\eta_n} \eta_n^{c_1 - 1/2} \int_{\ell=0}^{\lambda(\log n)} e^{-[\xi(\ell, \log n) - \xi(\ell_o, \log n)] r_n} [\xi(\ell, \log n) \xi(\ell_o, \log n)]^{1/2} d\ell =$

$c_2$,

(b) $\lim_{n\to\infty} e^{-\eta_n} \eta_n^{c_1-1/2} \sup_\ell \int_{t=\ell}^{\ell+r_n} e^{-[\xi(t,\log n)-\xi(\ell_o,\log n)]r_n} [\xi(t,\log n)\,\xi(\ell_o,\log n)]^{1/2}\, dt = 0$.

ASSUMPTION A6. $\lim_{u\to\infty} \inf_{\mathbf{x}:U(\mathbf{x})>u} \langle \mathbf{x}/|\mathbf{x}|, \nabla U(\mathbf{x})/|\nabla U(\mathbf{x})| \rangle > 0$.

In the following theorems let $N_n$ be a positive integer-valued random variable such that $N_n/n \xrightarrow{p} 1$. Also define

$$\mu_n = \frac{\eta_n - c_1 \log \eta_n - \log c_2^{-1}\sqrt{2\pi}}{\xi_n} \quad \text{and} \quad \sigma_n = 1/\xi_n, \tag{27}$$

and denote by

$$\Lambda(x) = \exp\{-e^{-x}\}, \quad -\infty < x < \infty,$$

the Gumbel distribution function.

THEOREM 4. *Assume that for each $\tau > 0$, Assumptions A1–A5 hold with $r_n$ defined by (26). Then*

$$(M_{\mathrm{NNG}(N_n)} - \mu_n)/\sigma_n \xrightarrow{d} \Lambda.$$

THEOREM 5. *Assume that for each $\tau > 0$, Assumptions A1–A6 hold with $r_n$ defined by (26). Then*

$$(M_{\mathrm{MST}(N_n)} - \mu_n)/\sigma_n \xrightarrow{d} \Lambda.$$

Clearly, if the conclusions of the theorems hold, then by weak convergence one has some flexibility in choosing the normalization; indeed, any $\tilde{\eta}_n$ and $\tilde{\xi}_n$ with

$$\tilde{\eta}_n - \eta_n \to 0 \quad \text{and} \quad \tilde{\xi}_n/\xi_n \to 1$$

can replace $\eta_n$ and $\xi_n$.

Assumptions A1–A6 are designed to meet not only analytic but also geometric considerations in the proofs, but still are quite general. Assumption A1 means that the normalized gradient of $U$ on a level curve $U(\mathbf{x}) = u$, where $u$ is large, will change gradually when the location changes gradually relative to the total length of the level curve. Assumptions A2 and A6 imply that the level curves are not highly asymmetrically shaped. Assumption A3 says that both $\lambda(u)$ and $|\nabla(u)|$ are $O$-regularly varying functions [cf. Bingham, Goldie and Teugels (1987)]. Assumption A4 is a very weak technical condition on the smoothness of $U$. Assumption A5 is the most significant condition in this group, which is aimed at connecting the normalizations for the longest edges on the graphs to the second-order Taylor expansion of $\xi(\ell, \log n)$ in areas where $\xi(\ell, \log n)$ is close to its minimum value $\xi(\ell_o, \log n)$. To further illustrate what these conditions mean and how to verify them, we proceed to examples in the next two sections.

**5. Homogeneous exponents.** The first general example considers distributions with exponents which are homogeneous in the sense that $f(\mathbf{x}) = e^{-U(\mathbf{x})} = e^{-V(\mathbf{x})-c}$ with $V(k\mathbf{x}) = k^\alpha V(\mathbf{x})$, where the level curves for $V$ are convex. In this case it is far more natural to use polar coordinates than the time-varying, locally orthogonal coordinate system in the general theory. In polar coordinates, $\mathbf{x} = (r\cos\theta, r\sin\theta)$; we may then write $V(\mathbf{x}) = r^\alpha g(\theta)$, so that

$$f(\mathbf{x}) = e^{-U(\mathbf{x})} = e^{-r^\alpha g(\theta)-c}. \tag{28}$$

Straightforward computations show that

$$\nabla U(\mathbf{x}) = r^{\alpha-1} A(\theta) \begin{bmatrix} \cos\theta \\ \sin\theta \end{bmatrix},$$

$$|\nabla U(\mathbf{x})| = r^{\alpha-1}|A(\theta)|^{1/2} =: r^{\alpha-1}h(\theta), \tag{29}$$

where

$$A(\theta) = \begin{bmatrix} \alpha g(\theta) & -g'(\theta) \\ g'(\theta) & \alpha g(\theta) \end{bmatrix} \quad \text{and} \quad |A(\theta)| = \alpha^2 g(\theta)^2 + g'(\theta)^2.$$

Defining $v_c(n) = \log n - c$, we have that $U(\mathbf{x}) = \log n$ is equivalent to $V(\mathbf{x}) = v_c(n)$. Hence, writing $r(\theta, \log n)$ for the solution to $U(\mathbf{x}) = r^\alpha g(\theta) + c = \log n$ and $\mathbf{x}(\theta, \log n)$ for the corresponding $\mathbf{x}$, we have that $r(\theta, \log n) = v_c(n)^{1/\alpha} g(\theta)^{-1/\alpha}$, and

$$\tilde{\xi}(\theta, \log n) = |\nabla U(\mathbf{x}(\theta, \log n))|$$
$$= v_c(n)^{1-1/\alpha} g(\theta)^{-1+1/\alpha} h(\theta) =: v_c(n)^{1-1/\alpha} k(\theta). \tag{30}$$

The final ingredient needed to treat distributions of the form (28) is to note that $\frac{\partial \ell(\theta)}{\partial \theta} = |\frac{\partial \mathbf{x}(\theta)}{\partial \theta}|$ and use

$$x(\theta, \log n) = r(\theta, \log n)(\cos\theta, \sin\theta) = v_c(n)^{1/\alpha} g(\theta)^{-1/\alpha}(\cos\theta, \sin\theta)$$

to obtain

$$\frac{\partial \ell(\theta)}{\partial \theta} = v_c(n)^{1/\alpha} g(\theta)^{-1/\alpha} \left[\left(\frac{g'(\theta)}{\alpha g(\theta)}\right)^2 + 1\right]^{1/2} =: v_c(n)^{1/\alpha} m(\theta). \tag{31}$$

THEOREM 6. *Assume that the density is of the form* (28) *with $\alpha > 1$ and $g(\theta)$ bounded away from zero and three times continuously differentiable on the torus $[0, 2\pi]$. Suppose further that $k(\theta)$ assumes its minimum value at $d$ distinct points $\theta_0, \ldots, \theta_{d-1}$. Then Assumptions A1–A6 hold with $\xi_n/\tilde{\xi}_n \to 1$ for $\tilde{\xi}_n = (\log n)^{1-1/\alpha} k(\theta_0)$, $\eta_n = \log_2 n$, $c_1 = 1$ and*

$$c_2 = \sqrt{2\pi} k(\theta_0)^{3/2} \sum_{i=0}^{d-1} k''(\theta_i)^{-1/2} m(\theta_i).$$

PROOF. If follows from (29) that

(32) $$\frac{\nabla U(\mathbf{x})}{|\nabla U(\mathbf{x})|} = \frac{A(\theta)}{|A(\theta)|^{1/2}} \begin{bmatrix} \cos\theta \\ \sin\theta \end{bmatrix}.$$

It is clear that Assumption A1 is equivalent to the continuity of the function on the right-hand side of (32), which follows from the assumptions. That Assumptions A2 and A3 hold are immediate consequences of (29). It may be seen that $U^{(i,j)} = V^{(i,j)} = r^{\alpha-(i+j)}g_{i,j}(\theta)$ for $i+j=2$ and suitable functions $g_{i,j}(\theta)$. Assumption A4 follows from this and (30). To show Assumption A6, note that by (32) we have

$$\left\langle \frac{\mathbf{x}}{|\mathbf{x}|}, \frac{\nabla U(\mathbf{x})}{|\nabla U(\mathbf{x})|} \right\rangle = \frac{\alpha g(\theta)}{|A(\theta)|^{1/2}},$$

which is bounded away from zero since $g$ is.

We now verify Assumption A5. Using in order (31), (30) and Erdelyi [(1956), page 37],

$$I_n := \int_\ell e^{-[\xi(\ell,\log n)-\xi(\ell_o,\log n)]r_n}[\xi(\ell,\log n)\xi(\ell_o,\log n)]^{1/2}\,d\ell$$

$$= \int_0^{2\pi} e^{-[\tilde{\xi}(\theta,\log n)-\tilde{\xi}(\theta_0,\log n)]r_n}[\tilde{\xi}(\theta,\log n)\tilde{\xi}(\theta_0,\log n)]^{1/2}v_c(n)^{1/\alpha}m(\theta)\,d\theta$$

$$= v_c(n)\int_0^{2\pi} e^{-r_n v_c(n)^{1-1/\alpha}[k(\theta)-k(\theta_0)]}[k(\theta)k(\theta_0)]^{1/2}m(\theta)\,d\theta$$

$$\sim v_c(n)k(\theta_0)\sqrt{2\pi}\sum_{i=0}^{d-1}(r_n v_c(n)^{1-1/\alpha}k''(\theta_i))^{-1/2}m(\theta_i).$$

Clearly $v_c(n) \sim \log n$ and by (25), (26) and (30),

$$r_n v_c(n)^{1-1/\alpha}k(\theta_0) = r_n \xi(\ell_0,\log n) \sim \eta_n.$$

Hence

$$I_n \sim \sqrt{2\pi}(\log n)k(\theta_0)^{3/2}\eta_n^{-1/2}\sum_{i=0}^{d-1}(k''(\theta_i))^{-1/2}m(\theta_i)$$

and Assumption A5(a) follows. Similar considerations prove Assumption A5(b). □

We next apply the result to classes of elliptically contoured distributions and to Weibull-type distributions.

EXAMPLE 1. Consider elliptically contoured distributions with log density

$$U(\mathbf{x}) = r^\alpha((\cos\theta)^2 - 2\rho\cos\theta\sin\theta + (\sin\theta)^2)/d + c.$$

In particular, this includes the bivariate normal with standardized marginals and correlation $\rho \in (-1, 1)$, but $\rho \neq 0$, which is obtained for $\alpha = 2$, $c_2 = -\log(2\pi(1-\rho^2)^{1/2})$, and $d = 2(1-\rho^2)$. These distributions are of the form (28) with

$$g(\theta) = ((\cos\theta)^2 - 2\rho\cos\theta\sin\theta + (\sin\theta)^2)/d = (1 - 2\rho\cos\theta\sin\theta)/d,$$

and hence with

$$k(\theta) = (1 - 2\rho\cos\theta\sin\theta)^{-1+1/\alpha}$$
$$\times \sqrt{\alpha^2(1-2\rho\cos\theta\sin\theta)^2 + (2\rho - 4\rho\cos^2\theta)^2}/d^{1/\alpha}.$$

The conditions of Theorem 6 are satisfied, and hence the results of Theorems 4 and 5 hold with

$$(33) \quad \mu_n = \frac{\log_2 n - \log_3 n - \log(c_2^{-1}\sqrt{2\pi})}{(\log n)^{1-1/\alpha} k(\theta_0)} \quad \text{and} \quad \sigma_n = \frac{1}{(\log n)^{1-1/\alpha} k(\theta_0)}.$$

Computation of $k(\theta_0)$ and $c_2$ requires computer algebra, with a slight simplification obtained by noting that the constant $c_2$ is independent of the value of $d$. We consider the normal case, which has $\alpha = 2$ and $d = 2(1-\rho^2)$, and present the results for a few values of $\rho$ in Table 1. By symmetry, the values for $-\rho$ are the same as for $\rho$.

EXAMPLE 2. In this example we consider independent marginals with Weibull-type densities $\text{const} \times e^{-|x_1|^\alpha - |x_2|^\alpha}$. The log density then is

$$U(\mathbf{x}) = r^\alpha(|\cos\theta|^\alpha + |\sin\theta|^\alpha) + c$$

and is of the form (28) with

$$g(\theta) = (|\cos\theta|^\alpha + |\sin\theta|^\alpha).$$

Hence,

$$k(\theta) = \alpha(|\cos\theta|^\alpha + |\sin\theta|^\alpha)^{1-1/\alpha}$$
$$\times \{(|\cos\theta|^\alpha + |\sin\theta|^\alpha)^2$$
$$\times (-\text{sign}(\cos\theta)\sin\theta|\cos\theta|^{\alpha-1} + \text{sign}(\sin\theta)\cos\theta|\sin\theta|^{\alpha-1})^2\}^{1/2}.$$

To assure three times differentiability we assume that $\alpha > 4$.

TABLE 1

| $\rho$ | 0.1 | 0.3 | 0.5 | 0.7 | 0.9 |
|---|---|---|---|---|---|
| $k(\theta_0)$ | 1.348837 | 1.240718 | 1.153867 | 1.084742 | 1.026251 |
| $c_2$ | 19.2383720 | 7.9460116 | 4.0933240 | 1.9501568 | 0.5414317 |

TABLE 2

| $\alpha$ | 5 | 6 | 7 | 8 |
|---|---|---|---|---|
| $k(\theta_0)$ | 0.769 | 0.472 | 0.280 | 0.162 |
| $c_2$ | 0.631 | 0.307 | 0.151 | 0.072 |

Again the conditions of Theorem 6 are satisfied, and hence the results of Theorems 4 and 5 hold with $\mu_n, \sigma_n$ given by (33) in the previous example. A few examples of values of $k(\theta_0)$, $c_2$, obtained with Maple, are given in Table 2.

It might be noted that even if the result requires $\alpha > 4$, it can be shown to hold also for $2 < \alpha \leq 4$.

**6. Parallel level curves.** Consider the situation where the level curves are parallel, namely they are given by

$$(34) \quad \{\mathbf{x} : U(\mathbf{x}) = u\} = \{\mathbf{c}(t) + \omega(u)\mathbf{n}(t) : t \in [0, L)\}, \quad u \geq \text{ some } u_o > 0,$$

where $L$ is a finite positive constant, $\mathbf{c}(t)$ is a closed, strictly convex curve parameterized clockwise by length, $\omega(u)$ is an increasing function with $\omega(u_o) = 0$ and $\omega(\infty) = \infty$, and $\mathbf{n}(t)$ is the unit normal of $\mathbf{c}$ at the parameter $t$. Assume that $\omega(u)$ is differentiable and $\mathbf{c}(t)$ is twice continuously differentiable with $|\ddot{\mathbf{c}}(t)| \in (0, \infty)$, where "·" refers to differentiation with respect to $t$. Since $|\dot{\mathbf{c}}(t)| \equiv 1$, it is easy to see that

$$(35) \quad \dot{\mathbf{n}}(t) = |\ddot{\mathbf{c}}(t)|\dot{\mathbf{c}}(t).$$

Also it is clear that

$$(36) \quad \xi(\ell, u) = |\nabla U(\mathbf{x}(\ell, u))| = \frac{1}{\omega'(u)} \quad \text{for all } \ell,$$

so that the asymptotic results here are different in nature from what was considered in Theorem 6. Since $\xi(\ell, u)$ does not depend on $\ell$ we will denote it by $\xi(u)$ henceforth. Define $t(\ell)$ to be the inverse function of

$$\ell(t) = \text{arc length of the level curve } (u = \log n)$$
$$\text{from } \mathbf{c}(0) + \omega(\log n)\mathbf{n}(0) \text{ to } \mathbf{c}(t) + \omega(\log n)\mathbf{n}(t)$$
$$= \int_0^t |\dot{\mathbf{c}}(v) + \omega(\log n)\dot{\mathbf{n}}(v)|\,dv = \int_0^t [1 + \omega(\log n)|\ddot{\mathbf{c}}(v)|]\,dv,$$

where the rightmost equality follows from (35) in conjunction with $|\dot{\mathbf{c}}(v)| \equiv 1$. Hence,

$$\mathbf{x}^{(1,0)}(\ell, u) = [\dot{\mathbf{c}}(t(\ell)) + \omega(u)\dot{\mathbf{n}}(t(\ell))]t'(\ell)$$

$$
\begin{align}
(37) \qquad &= [1 + \omega(u)|\ddot{\mathbf{c}}(t(\ell))|]t'(\ell)\dot{\mathbf{c}}(t(\ell)) \\
&= \frac{1 + \omega(u)|\ddot{\mathbf{c}}(t(\ell))|}{1 + \omega(\log n)|\ddot{\mathbf{c}}(t(\ell))|}\dot{\mathbf{c}}(t(\ell)).
\end{align}
$$

Also

$$
\begin{align}
\lambda(u) &= \int_0^L [1 + \omega(u)|\ddot{\mathbf{c}}(t)|]\, dt \\
(38) \qquad &= L + \omega(u)\int_0^L |\ddot{\mathbf{c}}(t)|\, dt \sim \omega(u)\int_0^L |\ddot{\mathbf{c}}(t)|\, dt.
\end{align}
$$

THEOREM 7. *Assume that* (34) *holds, where* $\mathbf{c}(t)$ *is twice continuously differentiable, and* $\omega(u) = u^\alpha \exp(\int_{y_o}^u \frac{a(y)}{y}\, dy)$ *for some* $\alpha > 0$, $a(y) \to 0$ *and* $ya'(y) \to 0$ *as* $y \to \infty$. *Then Assumptions* A1–A6 *hold with* $\eta_n = \log[\xi(\log n) \times \lambda(\log n)]$, $c_1 = 1/2$ *and* $c_2 = 1$.

PROOF. Since $\nabla U(\mathbf{x})/|\nabla U(\mathbf{x})| = \mathbf{n}(t)$, it is not difficult to see that Assumption A1 follows from the continuity of $\mathbf{n}(t)$. Both Assumptions A2 and A3 hold trivially in view of (36) and the assumption on $\omega$. We next verify Assumption A5. By (36) again,

$$
e^{-\eta_n}\eta_n^{c_1-1/2}\int_{\ell=0}^{\lambda(\log n)} e^{-[\xi(\ell,\log n)-\xi(\ell_o,\log n)]r_n}[\xi(\ell,\log n)\xi(\ell_o,\log n)]^{1/2}\, d\ell
$$
$$
= e^{-\eta_n}\eta_n^{c_1-1/2}\lambda(\log n)\xi(\log n).
$$

Hence, Assumption A5(a) is satisfied for the choice of constants in this theorem. The verification of Assumption A5(b) is entirely similar and therefore omitted. To verify Assumption A6, note that

$$
\left\langle \frac{\mathbf{x}}{|\mathbf{x}|}, \frac{\nabla U(\mathbf{x})}{|\nabla U(\mathbf{x})|} \right\rangle = \left\langle \frac{\mathbf{c}(t) + \omega(u)\mathbf{n}(t)}{|\mathbf{c}(t) + \omega(u)\mathbf{n}(t)|}, \mathbf{n}(t) \right\rangle \sim |\mathbf{n}(t)|^2 = 1 \qquad \text{as } u \to \infty.
$$

Hence Assumption A6 holds.

We finally verify Assumption A4. Solving $(\frac{\partial u}{\partial x_1}, \frac{\partial u}{\partial x_2})$ and $(\frac{\partial t}{\partial x_1}, \frac{\partial t}{\partial x_2})$ in the following in terms of $x_1$, $x_2$:

$$
1 = \frac{\partial x_1}{\partial x_1} = (\dot{c}_1 + \omega\dot{n}_1)\frac{\partial t}{\partial x_1} + \omega'n_1\frac{\partial u}{\partial x_1},
$$
$$
0 = \frac{\partial x_1}{\partial x_2} = (\dot{c}_1 + \omega\dot{n}_1)\frac{\partial t}{\partial x_2} + \omega'n_1\frac{\partial u}{\partial x_2},
$$
$$
0 = \frac{\partial x_2}{\partial x_1} = (\dot{c}_2 + \omega\dot{n}_2)\frac{\partial t}{\partial x_1} + \omega'n_2\frac{\partial u}{\partial x_1},
$$
$$
1 = \frac{\partial x_2}{\partial x_2} = (\dot{c}_2 + \omega\dot{n}_2)\frac{\partial t}{\partial x_2} + \omega'n_2\frac{\partial u}{\partial x_2},
$$

we get

$$\left(\frac{\partial u}{\partial x_1}, \frac{\partial u}{\partial x_2}\right) = \frac{1}{\omega'}(-\dot{c}_2, \dot{c}_1) \quad \text{and} \quad \left(\frac{\partial t}{\partial x_1}, \frac{\partial t}{\partial x_2}\right) = \frac{1}{1+\omega|\ddot{\mathbf{c}}|}(\dot{c}_1, \dot{c}_2).$$

Hence

$$\frac{\partial^2 u}{\partial x_1^2} = -\frac{\omega^{(2)}\dot{c}_2(\partial u/\partial x_1) - \omega'\ddot{c}_2(\partial t/\partial x_1)}{(\omega')^2} = \frac{\omega^{(2)}\dot{c}_2^2}{(\omega')^3} + \frac{\dot{c}_1\ddot{c}_2}{\omega'(1+\omega|\ddot{\mathbf{c}}|)},$$

$$\frac{\partial^2 u}{\partial x_2^2} = \frac{\omega^{(2)}\dot{c}_1(\partial u/\partial x_2) - \omega'\ddot{c}_1(\partial t/\partial x_2)}{(\omega')^2} = \frac{\omega^{(2)}\dot{c}_1^2}{(\omega')^3} + \frac{\ddot{c}_1\dot{c}_2}{\omega'(1+\omega|\ddot{\mathbf{c}}|)},$$

$$\frac{\partial^2 u}{\partial x_1 \partial x_2} = \frac{\omega^{(2)}\dot{c}_1(\partial u/\partial x_1) - \omega'\ddot{c}_1(\partial t/\partial x_1)}{(\omega')^2} = -\frac{\omega^{(2)}\dot{c}_1\dot{c}_2}{(\omega')^3} - \frac{\dot{c}_1\ddot{c}_1}{\omega'(1+\omega|\ddot{\mathbf{c}}|)}.$$

To show Assumption A4, by the computations above and the fact that $|\dot{\mathbf{c}}|$ and $|\ddot{\mathbf{c}}|$ are bounded, it suffices to show that

$$(39) \qquad \lim_{u\to\infty}\left(\frac{|\omega^{(2)}(u)|}{[\omega'(u)]^3} + \frac{1}{\omega(u)\omega'(u)}\right)(\omega'(u)\log u)^2 = 0.$$

The assumption of the theorem implies that both $\omega$ and $\omega'$ are regularly varying and in fact

$$\omega'(u) = (\alpha + a(y_o))u^{\alpha-1}\exp\left(\int_{y_o}^{u}\frac{\tilde{a}(y)}{y}\,dy\right),$$

where

$$\tilde{a}(y) = a(y) + \frac{ya'(y)}{\alpha + a(y)} \to 0 \quad \text{as } y \to \infty.$$

It is easily seen that

$$\frac{\omega'(u)}{\omega(u)} = \frac{\alpha + a(u)}{u} \quad \text{and} \quad \frac{\omega^{(2)}(u)}{\omega'(u)} = \frac{\alpha - 1 + \tilde{a}(u)}{u},$$

from which (39) is straightforward. □

EXAMPLE 3. An example of such a distribution is the bivariate normal distribution with independent standard normal marginals, in which case we can take $\mathbf{c}(t) = \mathbf{n}(t) = (\cos t, \sin t)$, $t \in [0, 2\pi)$, and

$$\omega(u) = \sqrt{2[u - \log(2\pi)]} - 1, \qquad u \geq \log(2\pi) + 1/2.$$

The conditions of Theorem 7 are satisfied and so Assumptions A1–A6 hold with

$$\xi_n = (2\log n)^{1/2}, \qquad \eta_n = \log 4\pi + \log_2 n, \qquad c_1 = 1/2, \qquad c_2 = 1.$$

Hence the results of Theorems 4 and 5 hold with

$$\mu_n = \frac{\log_2 n - (1/2)\log_3 n + \log(2\sqrt{2\pi})}{(2\log n)^{1/2}} \quad \text{and} \quad \sigma_n = \frac{1}{(2\log n)^{1/2}}.$$

This is consistent with the normalization obtained in Penrose (1998).

**7. Technical lemmas.** In this section we present three technical lemmas which are essential for Theorems 4 and 5. No significant continuity will be lost if a reader postpones the details in the proofs of these lemmas during the first reading.

The first one, Lemma 8, deduces three consequences of Assumption A4 which are more directly usable in the proofs of the theorems. Lemma 10 takes one step toward evaluating integrals like (1) and the last one, Lemma 9, approximates $F(S(\mathbf{x}; r_n))$, and in particular shows that appropriate sectors of $S(\mathbf{x}; r_n)$ may be neglected asymptotically.

In addition to previous notation, we will use in the proofs

$$(40) \quad u_{n,-b} = \log n - b \log_2 n, \quad u_{n,b} = \log n + b \log_2 n, \quad b > 0.$$

LEMMA 8. *Assume that Assumption A4 holds, namely $\lim_{U(\mathbf{x}) \to \infty}(|U^{(i,j)}(\mathbf{x})|/|\nabla U(\mathbf{x})|^2)(\log U(\mathbf{x}))^2 = 0$ for any $i, j \geq 0$ with $i + j = 2$, and consider the coordinate system based on the level curve $U(\mathbf{x}) = \log n$. Then, for any finite constant $b > 0$,*

$$(A4a) \quad \lim_{n \to \infty} \sup_{\mathbf{x}\,:\,|U(\mathbf{x}) - \log n| \leq b \log_2 n} \frac{|U^{(i,j)}(\mathbf{x})|}{|\nabla U(\mathbf{x})|^2} (\log_2 n)^2 = 0$$

*for any $i, j$ with $i + j = 2$.*

$$(A4b) \quad \lim_{n \to \infty} \sup_{|v - \log n| \leq b \log_2 n} \left|\frac{\xi(\ell, v)}{\xi(\ell, \log n)} - 1\right| \log_2 n = 0.$$

$$(A4c) \quad \lim_{n \to \infty} \sup_{|v - \log n| \leq b \log_2 n} ||\mathbf{x}^{(1,0)}(\ell, v)| - 1| = 0.$$

PROOF. Equation (A4a) is an immediate consequence of Assumption A4 since $\log U(\mathbf{x})/\log_2 n \to 1$ in the indicated range.

The proof of parts (A4b) and (A4c) uses the relation (18), that $\frac{\partial \mathbf{x}}{\partial u} = \nabla U(\mathbf{x})/|\nabla U(\mathbf{x})|^2$, in two ways. First, this relation implies that

$$(41) \quad \left|\frac{\partial x_1(\ell, v)}{\partial v}\right| \leq \frac{1}{|\nabla U(\mathbf{x})|} \quad \text{and} \quad \left|\frac{\partial x_2(\ell, v)}{\partial v}\right| \leq \frac{1}{|\nabla U(\mathbf{x})|},$$

and second, it is equivalent to

$$(42) \quad \mathbf{x}(\ell, v) - \mathbf{x}(\ell, \log n) = \int_{\log n}^{v} \frac{\nabla U(\mathbf{x}(\ell, s))}{|\nabla U(\mathbf{x}(\ell, s))|^2}\, ds.$$

Let
$$M_n = \sup_{|v-\log n|\leq b\log_2 n, i+j=2} \frac{|U^{(i,j)}(\mathbf{x}(\ell,v))|}{|\nabla U(\mathbf{x}(\ell,v))|^2},$$

so that $M_n(\log_2 n)^2 \to 0$ by (A4a). Now, to prove (A4b), note that, by straightforward differentiation,

$$\left|\frac{\partial \xi(\ell,v)}{\partial v}\right| = \left|\frac{\partial |\nabla U(\mathbf{x}(\ell,v))|}{\partial v}\right|$$
$$= \left|\frac{\partial}{\partial x_1}|\nabla U(\mathbf{x}(\ell,v))|\frac{\partial x_1}{\partial v} + \frac{\partial}{\partial x_2}|\nabla U(\mathbf{x}(\ell,v))|\frac{\partial x_2}{\partial v}\right|$$
$$= \frac{|\langle \nabla U, (U^{(2,0)}, U^{(1,1)})\rangle(\partial x_1/\partial v) + \langle \nabla U, (U^{(1,1)}, U^{(0,2)})\rangle(\partial x_2/\partial v)|}{|\nabla U|}.$$

Using (41), for $|v - \log n| \leq b\log_2 n$ we hence have that $|\frac{\partial \log \xi(\ell,v)}{\partial v}| \leq 4M_n$. Thus, by integration,

$$(e^{-4M_n b \log_2 n} - 1)\log_2 n \leq \left|\frac{\xi(\ell,v)}{\xi(\ell,\log n)} - 1\right|\log_2 n \leq (e^{4M_n b \log_2 n} - 1)\log_2 n$$

and (A4b) follows, since $M_n(\log_2 n)^2 \to 0$.

Next, to prove (A4c), we note that similar calculations as for (43) give that

$$\left|\frac{\partial}{\partial \ell}\frac{\nabla U(\mathbf{x}(\ell,v))}{|\nabla U(\mathbf{x}(\ell,v))|^2}\right| \leq 5M_n|\mathbf{x}^{(1,0)}(\ell,v)|.$$

Hence, interchanging differentiation and integration in (42), we obtain

$$|\mathbf{x}^{(1,0)}(\ell,v) - \mathbf{x}^{(1,0)}(\ell,\log n)| \leq 5M_n \int_{\log n}^{v} |\mathbf{x}^{(1,0)}(\ell,s)|\,ds.$$

By (19), $|\mathbf{x}^{(1,0)}(\ell,\log n)| = 1$, and hence

$$||\mathbf{x}^{(1,0)}(\ell,v)| - 1| \leq 5M_n \int_{\log n}^{v} (||\mathbf{x}^{(1,0)}(\ell,s)| - 1| + 1)\,ds.$$

It then follows from a Grönvall inequality that for $|v - \log n| \leq b\log_2 n$,

$$||\mathbf{x}^{(1,0)}(\ell,v)| - 1| \leq 5M_n b\log_2 n \exp\{5M_n b\log_2 n\} \to 0,$$

and (A4c) follows, since $M_n(\log_2 n) \to 0$. □

LEMMA 9. *Suppose that $U$ is differentiable with gradient $\nabla U$ and let, for some $\mathbf{x}$,*

(43) $$\varepsilon = \sup_{\mathbf{y}\in S(\mathbf{x};r)} |U(\mathbf{y}) - U(\mathbf{x}) - \langle \nabla U(\mathbf{x}), \mathbf{y} - \mathbf{x}\rangle|.$$

Let $\zeta \in (-1, 1]$. Then, with $\xi = |\nabla U(\mathbf{x})|$, there exist constants $\theta_1 = \theta_1(\mathbf{x}, r) \in [-1, 1]$ and $\frac{2}{1+\zeta} + \frac{1}{2} \leq \theta_2 = \theta_2(\mathbf{x}, r) \leq 0$ and such that

$$F(S(\mathbf{x}; r) \cap \{\mathbf{y} : \langle \mathbf{y} - \mathbf{x}, \nabla U(\mathbf{x})/|\nabla U(\mathbf{x})|\rangle \leq \zeta r\})$$
$$= (2\pi r)^{1/2} e^{-U(\mathbf{x})} \xi^{-3/2} e^{\xi r} e^{\theta_1 \varepsilon} [1 + \theta_2/(\xi r)].$$

PROOF. Let $B = \{\mathbf{y} : \langle \mathbf{y} - \mathbf{x}, \nabla U(\mathbf{x})/|\nabla U(\mathbf{x})|\rangle \leq \zeta r\}$. By (43), we can write

$$(44) \qquad F(S(\mathbf{x}; r) \cap B) = e^{-U(\mathbf{x})} e^{\theta_1 \varepsilon} \int_{S(\mathbf{x}; r) \cap B} e^{-\langle \nabla U(\mathbf{x}), \mathbf{y} - \mathbf{x}\rangle} \, d\mathbf{y}.$$

Then change variables with $\mathbf{y} - \mathbf{x} = A\mathbf{v}$ where

$$A = \begin{bmatrix} b & a \\ -a & b \end{bmatrix}$$

with $(a, b)' = \nabla U(\mathbf{x})/|\nabla U(\mathbf{x})|$. It is easy to see that

$$\int_{S(\mathbf{x}; r) \cap B} e^{-\langle \nabla U(\mathbf{x}), \mathbf{y} - \mathbf{x}\rangle} \, d\mathbf{y}$$
$$= \int_{S(0; r)} I(v_2 \leq \zeta r) e^{-\xi v_2} \, d\mathbf{v}$$
$$= \xi^{-2} \int_{S(0; \xi r)} I(v_2 \leq \xi \zeta r) e^{-v_2} \, d\mathbf{v}$$
$$= 2\xi^{-2} \int_{-\xi r}^{\xi \zeta r} e^{-v} [(\xi r)^2 - v^2]^{1/2} \, dv.$$

Next, letting $z = v + \xi r$, the previous expression is equal to

$$2\xi^{-2} e^{\xi r} \int_0^{\xi(1+\zeta)r} e^{-z} (2\xi r z - z^2)^{1/2} \, dz$$
$$= 2\xi^{-2} e^{\xi r} (2\xi r)^{1/2} \int_0^{\xi(1+\zeta)r} e^{-z} z^{1/2} (1 - z/(2\xi r))^{1/2} \, dz.$$

Since $1 - x/2 \leq (1-x)^{1/2} \leq 1$ for $x \in [0, 1]$, we have

$$\int_0^{\xi(1+\zeta)r} e^{-z} z^{1/2} \, dz - \frac{1}{4\xi r} \int_0^{\xi(1+\zeta)r} e^{-z} z^{3/2} \, dz$$
$$\leq \int_0^{\xi(1+\zeta)r} e^{-z} z^{1/2} (1 - z/(2\xi r))^{1/2} \, dz$$
$$\leq \int_0^{\infty} e^{-z} z^{1/2} \, dz = \Gamma(3/2).$$

The lower bound is

$$\Gamma(3/2) - \int_{\xi(1+\zeta)r}^{\infty} e^{-z} z^{1/2}\, dz - \frac{1}{4\xi r} \int_0^{\xi(1+\zeta)r} e^{-z} z^{3/2}\, dz,$$

and since $\Gamma(3/2) = \pi^{1/2}/2$, we obtain that

$$-\int_{\xi(1+\zeta)r}^{\infty} e^{-z} z^{1/2}\, dz - \frac{1}{4\xi r} \int_0^{\xi(1+\zeta)r} e^{-z} z^{3/2}\, dz$$

$$\leq \int_0^{\xi(1+\zeta)r} e^{-z} z^{1/2} (1 - z/(2\xi r))^{1/2}\, dz - \pi^{1/2}/2 \leq 0.$$

Now,

$$\int_{\xi(1+\zeta)r}^{\infty} e^{-z} z^{1/2}\, dz \leq \frac{1}{(\xi(1+\zeta)r)} \int_{\xi(1+\zeta)r}^{\infty} e^{-z/2} z^{3/2}\, dz \leq \frac{1}{\xi(1+\zeta)r} \Gamma(5/2)$$

and

$$\frac{1}{4\xi r} \int_0^{\xi(1+\zeta)r} e^{-z} z^{3/2}\, dz \leq \frac{1}{4\xi r} \Gamma(5/2).$$

Since $\Gamma(5/2) = 3\pi^{1/2}/4$ and $3\pi^{1/2}/8 \leq 1$, this concludes the proof. □

LEMMA 10. *Assume that Assumptions A2–A5 hold. Let constants $k \in [0, \infty)$ and $\delta_n \in [\underline{\delta}, \overline{\delta}] \subset (0, \infty)$. Then for any sufficiently large fixed $b \in (0, \infty)$,*

$$(45) \quad n\int_{\mathbf{x}} I(U(\mathbf{x}) \notin [u_{n,-b}, u_{n,b}]) e^{-\delta_n n F(S(\mathbf{x};r_n))} [nF(S(\mathbf{x};r_n))]^k\, dF(\mathbf{x}) \to 0,$$

*and uniformly for $\ell_{n,1} < \ell_{n,2} \in [0, \lambda(\log n)]$,*

$$n\int_{\mathbf{x}} I(U(\mathbf{x}) \in [u_{n,-b}, u_{n,b}], \ell(\mathbf{x}) \in [\ell_{n,1}, \ell_{n,2}])$$
$$(46) \qquad \times e^{-\delta_n n F(S(\mathbf{x};r_n))} [nF(S(\mathbf{x};r_n))]^k\, dF(\mathbf{x})$$
$$\sim \delta_n^{-1} \tau \frac{\int_{\ell=\ell_{n,1}}^{\ell_{n,2}} e^{-[\xi(\ell, \log n) - \xi(\ell_o, \log n)] r_n} \xi^{1/2}(\ell, \log n)\, d\ell}{\int_{\ell=0}^{\lambda(\log n)} e^{-[\xi(\ell, \log n) - \xi(\ell_o, \log n)] r_n} \xi^{1/2}(\ell, \log n)\, d\ell}.$$

*It follows from (45) and (46) that*

$$(47) \quad n\int_{\mathbf{x}} e^{-\delta_n n F(S(\mathbf{x};r_n))} [nF(S(\mathbf{x};r_n))]^k\, dF(\mathbf{x}) \sim \delta_n^{-1} \tau.$$

PROOF. In this proof we will assume for convenience of notation that $k = 0$ and $\delta_n \equiv 1$, since the extension to the general case is straightforward.

First note that by Assumption A2, (A4b) of Lemma 8 and the definition of $r_n$,

$$(48) \quad \sup_{u \in [u_{n,-b}, u_{n,b}], \ell} |\xi(\ell, u) - \xi(\ell, \log n)| r_n$$

$$\leq \sup_{u \in [u_{n,-b}, u_{n,b}], \ell} \left| \frac{\xi(\ell, u)}{\xi(\ell, \log n)} - 1 \right| O(\log_2 n) \to 0$$

and

$$\inf_{u \in [u_{n,-b}, u_{n,b}], \ell} \xi(\ell, u) r_n \geq \xi(\ell_o, \log n) r_n \to \infty, \qquad b > 0.$$

By Assumptions A4 and A2,

$$\lim_{n \to \infty} \sup_{\mathbf{x}: U(\mathbf{x}) \in [u_{n,-b}, u_{n,b}]} \sup_{\mathbf{y} \in S(\mathbf{x}; r_n)} |U(\mathbf{y}) - U(\mathbf{x}) - \langle \nabla U(\mathbf{x}), \mathbf{y} - \mathbf{x} \rangle| = 0, \qquad b > 0.$$

Hence, by Lemma 9 with $\zeta = 1$, (48) and (A4b) of Lemma 8, we have uniformly for all $\ell$ and $u \in [u_{n,-b}, u_{n,b}]$,

$$(49) \quad \begin{aligned} F(S(\mathbf{x}(\ell, u); r_n)) \\ \sim \sqrt{2\pi} e^{-u} \xi^{-3/2}(\ell, u) r_n^{1/2} e^{\xi(\ell, u) r_n} \\ \sim \sqrt{2\pi} e^{-u} \xi^{-3/2}(\ell, \log n) r_n^{1/2} e^{\xi(\ell_o, \log n) r_n} e^{[\xi(\ell, \log n) - \xi(\ell_o, \log n)] r_n} \\ = e^{-u} \chi_n(\ell), \end{aligned}$$

where

$$\chi_n(\ell) = c_2 \tau^{-1} \xi^{-3/2}(\ell, \log n) r_n^{1/2} e^{\eta_n} \eta_n^{-c_1} e^{[\xi(\ell, \log n) - \xi(\ell_o, \log n)] r_n}.$$

Now pick a large $u_o$ so that the bounds in Assumption A3 apply for $u > u_o$. By the nonintersection of level curves, Assumption A2 and (49), for all large $n$ we have, for some constant $b_1, b_2 > 0$,

$$\inf_{\mathbf{x}: u_o < U(\mathbf{x}) < u_{n,-b}} F(S(\mathbf{x}; r_n))$$

$$\geq \inf_\ell F(S(\mathbf{x}(\ell, u_{n,-b}); r_n))$$

$$\geq b_1 \inf_\ell e^{-u_{n,-b}} \xi^{-3/2}(\ell, \log n) r_n^{1/2} e^{\eta_n} \eta_n^{-c_1}$$

$$\geq b_2 n^{-1} (\log n)^b \xi^{-2}(\ell_o, \log n) e^{\eta_n} \eta_n^{1/2-c_1}.$$

It then follows from Assumption A3 that

$$\inf_{\mathbf{x}: u_o < U(\mathbf{x}) < u_{n,-b}} F(S(\mathbf{x}; r_n)) \geq b_3 n^{-1} (\log n)^{b_4}$$

for some constants $b_3, b_4$, where $b_4$ can be picked large provided that $b$ is large. Now for $\mathbf{x}$ for which $U(\mathbf{x}) \leq u_o$, it is clear that $F(S(\mathbf{x}; r_n))$ can be bounded below by the same bound, since the density there is bounded away from 0. Combining the two cases we conclude that if $b$ is large enough, we can pick $b_3$ and $b_4$ such that

$$(50) \quad n \int_{\mathbf{x}\,:\,U(\mathbf{x}) < u_{n,-b}} e^{-nF(S(\mathbf{x};r_n))} \, dF(\mathbf{x}) \leq n e^{-b_3 (\log n)^{b_4}} \to 0.$$

Next, by (20) and Assumption A3, there exists a constant $b_5$ such that

$$(51) \quad \begin{aligned} n \int_{\mathbf{x}\,:\,U(\mathbf{x}) > u_{n,b}} dF(\mathbf{x}) &= n \int_{u > u_{n,b}} \int_{\ell=0}^{\lambda(\log n)} e^{-u} \frac{|\mathbf{x}^{(1,0)}(\ell,u)|}{\xi(\ell,u)} \, d\ell \, du \\ &\leq n \int_{u > u_{n,b}} e^{-u} \frac{\lambda(u)}{\xi(\ell_o, u)} \, du \\ &\leq b_5 n \int_{u > u_{n,b}} e^{-u} u^{2\rho} \, du \\ &\sim b_5 n e^{-u_{n,b}} u_{n,b}^{2\rho} \to 0. \end{aligned}$$

Hence

$$(52) \quad n \int_{\mathbf{x}\,:\,U(\mathbf{x}) > u_{n,b}} e^{-nF(S(\mathbf{x};r_n))} \, dF(\mathbf{x}) \to 0$$

for a sufficiently large $b$. Now, (50) and (52) imply (45).

It follows from (24), (49), and (A4b) and (A4c) of Lemma 8 that

$$n \int I(U(\mathbf{x}) \in [u_{n,-b}, u_{n,b}], \ell(\mathbf{x}) \in [\ell_{n,1}, \ell_{n,2}]) e^{-nF(S(\mathbf{x};r))} \, dF(\mathbf{x})$$

$$\sim n \int_{u=u_{n,-b}}^{u_{n,b}} \int_{\ell=\ell_{n,1}}^{\ell_{n,2}} e^{-n(1+o(1))e^{-u}\chi_n(\ell)} e^{-u} \frac{1}{\xi(\ell, \log n)} \, d\ell \, du.$$

Make a change of variables in the above integral with

$$v = u - \log n - \log \chi_n(\ell).$$

By Assumptions A2, A3 and the assumption that $\eta_n = \log_2(n)$ in Assumption A5, we conclude that $\log \chi_n(\ell) = O(\log_2 n)$. As a result, it is possible to choose large enough $b$ such that $v_{n,-b} \to -\infty$ and $v_{n,b} \to \infty$, where $v_{n,\pm b}$ corresponds to $u_{n,\pm b}$ with respect to the above variable change. Since

$$\int_{-\infty}^{\infty} e^{-e^{-v}} e^{-v} \, dv = 1,$$

we have

$$n \int I(U(\mathbf{x}) \in [u_{n,-b}, u_{n,b}], \ell(\mathbf{x}) \in [\ell_{n,1}, \ell_{n,2}]) e^{-nF(S(\mathbf{x};r))} \, dF(\mathbf{x})$$

$$\sim \int_{v=v_{n,-b}}^{v_{n,b}} \int_{\ell=\ell_{n,1}}^{\ell_{n,2}} e^{-(1+o(1))e^{-v}} e^{-v} \frac{1}{\xi(\ell, \log n) \chi_n(\ell)} \, d\ell \, dv$$

$$\sim c_2^{-1} \tau r_n^{-1/2} e^{-\eta_n} \eta_n^{c_1} \int_{\ell=\ell_{n,1}}^{\ell_{n,2}} e^{-[\xi(\ell, \log n) - \xi(\ell_o, \log n)] r_n} \xi^{1/2}(\ell, \log n) \, d\ell$$

$$\sim c_2^{-1} \tau e^{-\eta_n} \eta_n^{c_1 - 1/2}$$

$$\times \int_{\ell=\ell_{n,1}}^{\ell_{n,2}} e^{-[\xi(\ell, \log n) - \xi(\ell_o, \log n)] r_n} [\xi(\ell, \log n) \xi(\ell_o, \log n)]^{1/2} \, d\ell$$

$$\sim \tau \frac{\int_{\ell=\ell_{n,1}}^{\ell_{n,2}} e^{-[\xi(\ell, \log n) - \xi(\ell_o, \log n)] r_n} \xi^{1/2}(\ell, \log n) \, d\ell}{\int_{\ell=0}^{\lambda(\log n)} e^{-[\xi(\ell, \log n) - \xi(\ell_o, \log n)] r_n} \xi^{1/2}(\ell, \log n) \, d\ell}$$

by Assumption A5. This proves (46) for $k = 0$ and $\delta_n = 1$. □

**8. Proofs of Theorems 4 and 5.** We continue to use below the notation of $u_{n,\pm b}$ defined in (40).

PROOF OF THEOREM 4. First consider the case $N_n \sim \text{Poisson}(n)$. We start by defining the sets $A_{n,1}, \ldots, A_{n,k_n}$ in Theorem 1. For convenience write

$$I_n(\ell) = \frac{\tau \int_{t=0}^{\ell} e^{-[\xi(t, \log n) - \xi(\ell_o, \log n)] r_n} \xi^{1/2}(t, \log n) \, dt}{\int_{t=0}^{\lambda(\log n)} e^{-[\xi(t, \log n) - \xi(\ell_o, \log n)] r_n} \xi^{1/2}(t, \log n) \, dt}, \qquad \ell \in [0, \lambda(\log n)).$$

Thus, by (46) of Lemma 10 with $k = 0$ and $\delta_n \equiv 1$, for some $b > 0$ we have uniformly for $0 \leq \ell_{n,1} < \ell_{n,2} < \lambda(\log n)$,

(53) $\quad \mu_n^{(1)}(\{\mathbf{x}(\ell, u) : \ell \in [\ell_{n,1}, \ell_{n,2}), u \in [u_{n,-b}, u_{n,b}]\}, r_n) \sim I_n(\ell_{n,2}) - I_n(\ell_{n,1})$

where $\mu_n^{(1)}$ is defined in (1). In the following we will continue to work with this choice of $b$. By Assumptions A5(a) and A5(b),

$$\varepsilon_n := \sup_{\ell} [I_n(\ell + r_n) - I_n(\ell)] \to 0.$$

Let $k_n$ be a sequence of integers such that

$$k_n \to \infty \quad \text{and} \quad k_n \varepsilon_n \to 0;$$

let $j_0 = 0$ and, for $i = 1, \ldots, k_n$, let $j_i$ be the largest positive integer $j$ such that $I_n(jr_n) \leq i\tau/k_n$, where $j_{k_n}$ is simply the largest positive integer $j$ such that $jr_n \leq \lambda(\log n)$. Note that since $\varepsilon_n = o(1/k_n)$, the differences between successive $j_i$'s tend to $\infty$ uniformly so that $j_{i-1} < j_i - 3$ for all large $n$. Define

$$A_{n,i} = \{\mathbf{x}(\ell, u) : \ell \in [j_{i-1} r_n, (j_i - 3) r_n), u \in [u_{n,-b}, u_{n,b}]\}, \qquad 1 \leq i \leq k_n,$$

and
$$A_{n,k_n+1} = A_{n,1}.$$

Also let
$$\bar{A}_{n,i} = \{\mathbf{x}(\ell, u) : \ell \in [(j_i - 3)r_n, j_i r_n), u \in [u_{n,-b}, u_{n,b}]\}, \qquad 1 \leq i \leq k_n - 1,$$

and
$$\bar{A}_{n,k_n} = \{\mathbf{x}(\ell, u) : \ell \in [(j_{k_n} - 3)r_n, \lambda(\log n)), u \in [u_{n,-b}, u_{n,b}]\}.$$

By the choice of the sets, it follows from (53) and the definition of $\varepsilon_n$ that

(54) $\quad 1/k_n - 4\varepsilon_n \leq \mu_n^{(1)}(A_{n,i}, r_n) \leq 1/k_n + \varepsilon_n, \qquad 1 \leq i \leq k_n,$

(55) $\quad \max_{1 \leq i \leq k_n - 1} \mu_n^{(1)}(\bar{A}_{n,i}, r_n) \leq 3\varepsilon_n \quad \text{and} \quad \mu_n^{(1)}(\bar{A}_{n,k_n}, r_n) \leq 3\varepsilon_n + 1/k_n,$

and hence that

(56) $\quad \max_{1 \leq i \leq k_n - 1} \frac{\mu_n^{(1)}(\bar{A}_{n,i}, r_n)}{\mu_n^{(1)}(A_{n,i}, r_n)} \leq \frac{3\varepsilon_n}{1/k_n - 4\varepsilon_n} \to 0.$

We now proceed to verify the conditions (a)–(e) of Theorem 1.

By convexity of the level curves,
$$c_{n,i} := \inf(|\mathbf{x} - \mathbf{y}| : \mathbf{x} \in A_{n,i}, \mathbf{y} \in A_{n,i+1}) = |\mathbf{x}((j_i - 3)r_n, u_{n,-b}) - \mathbf{x}(j_i r_n, u_{n,-b})|.$$

By (20) and the condition (A4c) of Lemma 8, the length $a_{n,i}$ of the arc that connects $\mathbf{x}((j_i - 3)r_n, u_{n,-b})$ with $\mathbf{x}(j_i r_n, u_{n,-b})$ on $\{\mathbf{x} : U(\mathbf{x}) = u_{n,-b}\}$ is asymptotically $3r_n$. By the same token, $\lambda(u_{n,-b}) \sim \lambda(\log n)$. Let $\theta_{n,i}$ be the angle between $\frac{\nabla U(\mathbf{x}((j_i-3)r_n, u_{n,-b}))}{|\nabla U(\mathbf{x}((j_i-3)r_n, u_{n,-b}))|}$ and $\frac{\nabla U(\mathbf{x}(j_i r_n, u_{n,-b}))}{|\nabla U(\mathbf{x}(j_i r_n, u_{n,-b}))|}$. Since $r_n = o(\lambda(\log n))$ by Assumptions A5(a) and A5(b), Assumption A1 then guarantees that $\max_i \theta_{n,i} \to 0$. Consider the triangle with base equal to the line that connects $\mathbf{x}((j_i - 3)r_n, u_{n,-b})$ and $\mathbf{x}(j_i r_n, u_{n,-b})$, and with sides determined by the tangents at these two points on the curve $U(\mathbf{x}) = u_{n,-b}$; let the two base angles be $\theta_1$ and $\theta_2$, say, and the corresponding two sides have lengths $s_1$, $s_2$, respectively. Figure 1 depicts what is described here.

Then
$$\frac{(s_1 + s_2)\cos\theta_{n,i}}{c_{n,i}} = \frac{(s_1 + s_2)\cos(\theta_1 + \theta_2)}{c_{n,i}} \leq \frac{s_1 \cos\theta_1 + s_2 \cos\theta_2}{c_{n,i}} = 1,$$

and so
$$\frac{a_{n,i}}{c_{n,i}} \leq \frac{s_1 + s_2}{c_{n,i}} \leq \frac{1}{\cos\theta_{n,i}} \to 1 \qquad \text{uniformly as } n \to \infty.$$

Consequently, $\min_i c_{n,i} > 2r_n$ for large $n$. Hence we have shown that adjacent $A_{n,i}$'s are at least $2r_n$ apart for large $n$. That nonadjacent $A_{n,i}$'s are at least $2r_n$ apart for large $n$ can also be established by the convexity of the level curves. Hence the condition (a) of Theorem 1 is proved.

It follows from (47) of Lemma 10 with $\delta_n \equiv 1$ and $k = 0$ that the condition (b) of Theorem 1 holds. To prove (c) of Theorem 1, note that it suffices to show that

$$\mu_n^{(1)}\{\mathbf{x}(\ell, u) : u \notin [u_{n,-b}, u_{n,b}]\} \to 0 \quad \text{and} \quad \sum_{j=1}^{k_n} \mu_n^{(1)}(\bar{A}_{n,i}, r_n) \to 0.$$

The first convergence follows from (45) of Lemma 10 with $\delta_n \equiv 1$ and $k = 0$, while the second convergence follows from (55) and (56), in conjunction with the condition (b) already proved above. Condition (d) of Theorem 1 holds by (54).

Finally we prove (e) of Theorem 1. By the simple inequality

$$P(A \cup B) \geq (P(A) + P(B))/2$$

we have

$$\mu_n^{(2)}(A, r_n) \leq \left(n \int_A e^{-nF(S(\mathbf{x};r_n))/2} \, dF(\mathbf{x})\right)^2$$

so that

$$\frac{\mu_n^{(2)}(A_{n,i}, r_n)}{\mu_n^{(1)}(A_{n,i}, r_n)} \leq \frac{(n \int_{A_{n,i}} e^{-nF(S(\mathbf{x};r_n))/2} \, dF(\mathbf{x}))^2}{n \int_{A_{n,i}} e^{-nF(S(\mathbf{x};r_n))} \, dF(\mathbf{x})}.$$

By (46) of Lemma 10 with $\delta_n \equiv 1/2, 1$ and $k = 0$, and (54),

$$\frac{(n \int_{A_{n,i}} e^{-nF(S(\mathbf{x};r_n))/2} \, dF(\mathbf{x}))^2}{n \int_{A_{n,i}} e^{-nF(S(\mathbf{x};r_n))} \, dF(\mathbf{x})} \sim 4\tau \mu_n^{(1)}(A_{n,i}, r_n) \to 0$$

uniformly in $i$, where $\tau$ is the constant specified in the theorem in defining $r_n$. Hence (e) of Theorem 1 is proved. Replacing $\tau$ by $e^{-x}$ completes the proof for the case $N_n \sim \text{Poisson}(n)$.

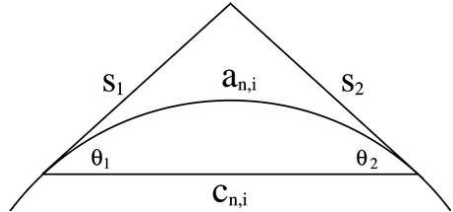

Fig. 1.

Next,
$$\begin{aligned}P(M_{\text{NNG}(N_n)} &\leq r_n) \\ &\leq P(M_{\text{NNG}(n)} \leq r_n) + P(M_{\text{NNG}(N_n)} \leq r_n < M_{\text{NNG}(n)}) \\ &\leq P(M_{\text{NNG}(N_n)} \leq r_n) + P(M_{\text{NNG}(n)} \leq r_n < M_{\text{NNG}(N_n)}) \\ &\quad + P(M_{\text{NNG}(N_n)} \leq r_n < M_{\text{NNG}(n)}).\end{aligned}$$

Using the inequality $1 - x \leq e^{-x}$, $x \leq 1$, it follows from (47) of Lemma 10 with $k = 0$, $\delta_n \equiv 1$ and $k = 1$, $\delta_n \equiv \delta$, where $\delta \in (0, 1)$, that (12) and (13) hold. Taking limits throughout the preceding inequalities, it follows from Theorem 2 and the first part of the proof with $N_n \sim \text{Poisson}(n)$ that
$$\lim_{n \to \infty} P(M_{\text{NNG}(n)} \leq r_n) = \lim_{n \to \infty} P(M_{\text{NNG}(N_n)} \leq r_n) = e^{-\tau}.$$

Applying the above argument again gives
$$\lim_{n \to \infty} P(M_{\text{NNG}(N_n)} \leq r_n) = \lim_{n \to \infty} P(M_{\text{NNG}(n)} \leq r_n) = e^{-\tau}$$

for any general $N_n$ satisfying $N_n/n \xrightarrow{p} 1$. Replacing $\tau$ by $e^{-x}$ completes the proof. □

PROOF OF THEOREM 5. It follows from Theorems 3 and 4 that
$$\lim_{n \to \infty} P(M_{\text{MST}(N_n)} \leq r_n) = \lim_{n \to \infty} P(M_{\text{NNG}(N_n)} \leq r_n) = e^{-\tau}$$

provided we show (16). To do that, it suffices to show that

(57) $$n \int \bar{F}^{\delta_n n}(S(\mathbf{x}; r_n)) I(|\mathbf{x}| > r_n/2) \, dF(\mathbf{x}) \to \tau,$$

and

(58) $$n \int \bar{F}^{\delta_n n}(S(\mathbf{x}; r_n) \cap S(0; |\mathbf{x}|)) I(|\mathbf{x}| > r_n/2) \, dF(\mathbf{x}) \to \tau$$

for any sequence $\delta_n \to 1$. For any fixed $\epsilon > 1$, it follows from Assumption A6 that for all large $n$,

(59) $$\begin{aligned}2|\mathbf{x}(\ell, u_{n,-b})| &\geq |\mathbf{x}(\ell, u_{n,-b})| + |\mathbf{x}(\ell, u_{n,-\epsilon b})| \\ &\geq |\mathbf{x}(\ell, u_{n,-b}) - \mathbf{x}(\ell, u_{n,-\epsilon b})|.\end{aligned}$$

By (18) and (A4b) of Lemma 8,
$$\begin{aligned}|\mathbf{x}(\ell, u_{n,-b}) - \mathbf{x}(\ell, u_{n,-\epsilon b})| &= \left|\int_{u_{n,-\epsilon b}}^{u_{n,-b}} \frac{\nabla U(\mathbf{x}(\ell, v))}{|\nabla U(\mathbf{x}(\ell, v))|^2} \, dv\right| \\ &\sim \frac{1}{\xi^2(\ell, \log n)} \left|\int_{u_{n,-\epsilon b}}^{u_{n,-b}} \nabla U(\mathbf{x}(\ell, v)) \, dv\right|.\end{aligned}$$

Now write
$$\int_{u_{n,-\epsilon b}}^{u_{n,-b}} \nabla U(\mathbf{x}(\ell, v)) \, dv = \int_{u_{n,-\epsilon b}}^{u_{n,-b}} \nabla U(\mathbf{x}(\ell, \log n)) \, dv$$
$$+ \int_{u_{n,-\epsilon b}}^{u_{n,-b}} [\nabla U(\mathbf{x}(\ell, v)) - \nabla U(\mathbf{x}(\ell, \log n))] \, dv.$$

Since
$$\left| \int_{u_{n,-\epsilon b}}^{u_{n,-b}} \nabla U(\mathbf{x}(\ell, \log n)) \, dv \right| = \xi(\ell, \log n)(\epsilon - 1) b \log_2 n,$$

if we can show

(60) $$\sup_{v \in [u_{n,-\epsilon b}, u_{n,-b}]} \frac{|\nabla U(\mathbf{x}(\ell, v)) - \nabla U(\mathbf{x}(\ell, \log n))|}{\xi(\ell, \log n)} \to 0,$$

then it follows at once that

(61) $$|\mathbf{x}(\ell, u_{n,-b}) - \mathbf{x}(\ell, u_{n,-\epsilon b})| \sim \frac{(\epsilon - 1) b \log_2 n}{\xi(\ell, \log n)}.$$

We now show (60). By the mean value theorem,

(62) $$\sup_{v \in [u_{n,-\epsilon b}, u_{n,-b}]} \frac{|\nabla U(\mathbf{x}(\ell, v)) - \nabla U(\mathbf{x}(\ell, \log n))|}{\xi(\ell, \log n)}$$
$$\leq \frac{\epsilon b \log_2 n}{\xi(\ell, \log n)} \left( \sup_{v \in [u_{n,-\epsilon b}, u_{n,-b}]} \left| \frac{\partial}{\partial v} U^{(1,0)}(\mathbf{x}(\ell, v)) \right| \right.$$
$$\left. + \sup_{v \in [u_{n,-\epsilon b}, u_{n,-b}]} \left| \frac{\partial}{\partial v} U^{(0,1)}(\mathbf{x}(\ell, v)) \right| \right).$$

As in the proof of (A4b) in Lemma 8, by (41) we have
$$\left| \frac{\partial}{\partial v} U^{(1,0)}(\mathbf{x}(\ell, v)) \right| = \left| U^{(2,0)}(\mathbf{x}(\ell, v)) \frac{\partial x_1(\ell, v)}{\partial v} + U^{(1,1)}(\mathbf{x}(\ell, v)) \frac{\partial x_2(\ell, v)}{\partial v} \right|$$
$$\leq \frac{1}{|\nabla U(\mathbf{x}(\ell, v))|} (|U^{(2,0)}(\mathbf{x}(\ell, v))| + |U^{(1,1)}(\mathbf{x}(\ell, v))|),$$

which, by (A4a) and (A4b), gives
$$\sup_{v \in [u_{n,-\epsilon b}, u_{n,-b}]} \left| \frac{\partial}{\partial v} U^{(1,0)}(\mathbf{x}(\ell, v)) \right| = o\left( \frac{\xi(\ell, \log n)}{(\log_2 n)^2} \right).$$

The same conclusion can be reached for $\sup_{v \in [u_{n,-\epsilon b}, u_{n,-b}]} |\frac{\partial}{\partial v} U^{(0,1)}(\mathbf{x}(\ell, v))|$. Thus, (60) follows from (62), and (61) is proved. Since $r_n \leq \eta_n / \xi_n$ for large $n$ where $\eta_n = O(\log_2 n)$, we conclude by (59), (60) and Assumption A2 that

$$\liminf_{n \to \infty} \inf_\ell |\mathbf{x}(\ell, u_{n,-b})| / r_n \geq cb$$

for some finite constant $c$ independent of $b$. Hence we can choose $b$ sufficiently large to ensure that $\inf_\ell |\mathbf{x}(\ell, u_{n,-b})|/r_n \geq 1$. By this and (45) of Lemma 10 it is the case that for a large enough $b$,

$$n \int \bar{F}^{\delta_n n}(S(\mathbf{x}; r_n)) I(|\mathbf{x}| \leq r_n/2) \, dF(\mathbf{x})$$

$$\leq n \int_{\mathbf{x}:\, U(\mathbf{x}) < u_{n,-b}} \bar{F}^{\delta_n n}(S(\mathbf{x}; r_n)) \, dF(\mathbf{x}) \to 0$$

so that (57) follows from Lemma 10.

We next consider the proof of (58). We will do so by proving that for a $b > 0$, sufficiently large,

(63) $\quad n \int_{\mathbf{x}:\, U(\mathbf{x}) > u_{n,b}} \bar{F}^{\delta_n n}(S(\mathbf{x}; r_n) \cap S(0; |\mathbf{x}|)) I(|\mathbf{x}| > r_n/2) \, dF(\mathbf{x}) \to 0,$

(64) $\quad n \int_{\mathbf{x}:\, U(\mathbf{x}) < u_{n,-b}} \bar{F}^{\delta_n n}(S(\mathbf{x}; r_n) \cap S(0; |\mathbf{x}|)) I(|\mathbf{x}| > r_n/2) \, dF(\mathbf{x}) \to 0$

and

$$n \int_{\mathbf{x}:\, u_{n,-b} \leq U(\mathbf{x}) \leq u_{n,b}} \bar{F}^{\delta_n n}(S(\mathbf{x}; r_n) \cap S(0; |\mathbf{x}|)) I(|\mathbf{x}| > r_n/2) \, dF(\mathbf{x}) \to \tau.$$
(65)

Clearly (63) follows from (51) in the proof of Lemma 10. To deal with (64), first write

$$n \int_{\mathbf{x}:\, U(\mathbf{x}) < u_{n,-b}} \bar{F}^{\delta_n n}(S(\mathbf{x}; r_n) \cap S(0; |\mathbf{x}|)) I(|\mathbf{x}| > r_n/2) \, dF(\mathbf{x})$$

$$= n \int_{\mathbf{x}:\, U(\mathbf{x}) < u_{n,-b}/2} \bar{F}^{\delta_n n}(S(\mathbf{x}; r_n) \cap S(0; |\mathbf{x}|)) I(|\mathbf{x}| > r_n/2) \, dF(\mathbf{x})$$

$$+ n \int_{\mathbf{x}:\, u_{n,-b}/2 \leq U(\mathbf{x}) < u_{n,-b}} \bar{F}^{\delta_n n}(S(\mathbf{x}; r_n) \cap S(0; |\mathbf{x}|)) I(|\mathbf{x}| > r_n/2) \, dF(\mathbf{x}).$$

We will show that both terms on the right-hand side tend to 0. Since $f(\mathbf{x}) = e^{-U(\mathbf{x})}$, by Assumptions A2, A3 and the mean value theorem, there exists some $b_1 \in (0, \infty)$ such that

$$\inf_{\mathbf{y} \in S(\mathbf{x}; r_n),\, U(\mathbf{x}) < u_{n,-b}/2} f(\mathbf{y}) \geq e^{-u_{n,-b}/2 - b_1 \xi(\ell_o, \log n) r_n} \geq e^{-u_{n,-b}/2 - b_1 \eta_n}.$$

Also observe that for $\mathbf{x}$ such that $|\mathbf{x}| > r_n/2$,

$$\text{area}(S(\mathbf{x}; r_n) \cap S(0; |\mathbf{x}|)) \geq b_2 \pi r_n^2$$

for some $b_2 \in (0, 1)$. Hence, for all large $n$,

$$\inf_{\mathbf{x}:\, |\mathbf{x}| > r_n/2,\, U(\mathbf{x}) < u_{n,-b}/2} F(S(\mathbf{x}; r_n) \cap S(0; |\mathbf{x}|)) \geq b_2 \pi r_n^2 e^{-u_{n,-b}/2 - b_1 \eta_n}$$

and hence

$$n \int_{\mathbf{x}:\, U(\mathbf{x}) < u_{n,-b}/2} \bar{F}^{\delta_n n}(S(\mathbf{x}; r_n) \cap S(0; |\mathbf{x}|)) I(|\mathbf{x}| > r_n/2) \, dF(\mathbf{x}) \tag{66}$$

$$\leq n \exp\{-\delta_n n b_2 \pi r_n^2 e^{-u_{n,-b}/2 - b_1 \eta_n}\} \to 0.$$

By Assumption A6, there exists some constant $\zeta \in (-1, 0)$ such that for all $\mathbf{x}$ with $U(\mathbf{x}) > u_{n,-b}/2$,

$$S(\mathbf{x}; r_n) \cap \left\{ \mathbf{y} : \left\langle \mathbf{y} - \mathbf{x}, \frac{\nabla U(\mathbf{x})}{|\nabla U(\mathbf{x})|} \right\rangle \leq \zeta r_n \right\} \subset S(\mathbf{x}; r_n) \cap S(0; |\mathbf{x}|).$$

Hence

$$\inf_{\mathbf{x}:\, u_{n,-b}/2 \leq U(\mathbf{x}) < u_{n,-b}} F(S(\mathbf{x}; r_n) \cap S(0; |\mathbf{x}|))$$

$$\tag{67} \geq \inf_{\mathbf{x}:\, u_{n,-b}/2 \leq U(\mathbf{x}) < u_{n,-b}} F\left(S(\mathbf{x}; r_n) \cap \left\{ \mathbf{y} : \left\langle \mathbf{y} - \mathbf{x}, \frac{\nabla U(\mathbf{x})}{|\nabla U(\mathbf{x})|} \right\rangle \leq \zeta r_n \right\}\right)$$

$$\geq \inf_{\ell} F\bigg( S(\mathbf{x}(\ell, u_{n,-b}); r_n)$$

$$\cap \left\{ \mathbf{y} : \left\langle \mathbf{y} - \mathbf{x}(\ell, u_{n,-b}), \frac{\nabla U(\mathbf{x}(\ell, u_{n,-b}))}{|\nabla U(\mathbf{x}(\ell, u_{n,-b}))|} \right\rangle \leq \zeta r_n \right\} \bigg).$$

It follows from Lemma 9 that uniformly for $u \in [u_{n,-b}, u_{n,b}]$ and $\ell$,

$$F\bigg( S(\mathbf{x}(\ell, u); r_n) \cap \left\{ \mathbf{x} : \left\langle \mathbf{x} - \mathbf{x}(\ell, u), \frac{\nabla U(\mathbf{x}(\ell, u))}{|\nabla U(\mathbf{x}(\ell, u))|} \right\rangle \leq \zeta r_n \right\} \bigg) \tag{68}$$

$$\sim F(S(\mathbf{x}(\ell, \log n); r_n)).$$

By (67) and (68), the same proof that leads to (50) in Lemma 10 now proves

$$n \int_{\mathbf{x}:\, u_{n,-b}/2 \leq U(\mathbf{x}) < u_{n,-b}} \bar{F}^{\delta_n n}(S(\mathbf{x}; r_n) \cap S(0; |\mathbf{x}|)) \tag{69}$$

$$\times I(r_n/2 < |\mathbf{x}| < u_{n,-b}) \, dF(\mathbf{x}) \to 0.$$

Hence (64) follows from (66) and (69). Making use of (68), the proof of (65) mirrors that of (46) of Lemma 10 and is omitted. □

**Acknowledgments.** Both authors wish to thank the referee whose most detailed report pointed out a number of flaws in the previous version of this paper and provided helpful ideas for the revision. Tailen Hsing also thanks Claudia Klüppelberg for her support of this research and her hospitality during his visit to Munich in 2002.

# REFERENCES

BINGHAM, N. H., GOLDIE, C. M. and TEUGELS, J. L. (1987). *Regular Variation.* Cambridge Univ. Press. MR898871

ERDELYI, A. (1956). *Asymptotic Expansions.* Dover, New York. MR78494

IBRAGIMOV, I. A. and LINNIK, YU. V. (1971). *Independent and Stationary Sequences of Random Variables.* Wolters-Noordhoff, Groningen. MR322926

KALLENBERG, O. (1983). *Random Measures.* Academic Press, New York. MR818219

KESTEN, H. and LEE, S. (1996). The central limit theorem for weighted minimal spanning trees on random points. *Ann. Appl. Probab.* **6** 495–527. MR1398055

LEADBETTER, M. R., LINDGREN, G. and ROOTZÉN, H. (1983). *Extremes and Related Properties of Random Sequences and Processes.* Springer, New York. MR691492

LEE, S. (1997). The central limit theorem for Euclidean minimal spanning trees. I. *Ann. Appl. Probab.* **7** 996–1020. MR1484795

PENROSE, M. D. (1997). The longest edge of the minimal spanning tree. *Ann. Appl. Probab.* **7** 340–361. MR1442317

PENROSE, M. D. (1998). Extremes for the minimal spanning tree on normally distributed points. *Adv. in Appl. Probab.* **30** 628–639. MR1663521

PENROSE, M. (2000). Central limit theorems for $k$-nearest neighbour distances. *Stochastic Process. Appl.* **85** 295–320. MR1731028

DEPARTMENT OF STATISTICS  
OHIO STATE UNIVERSITY  
COLUMBUS, OHIO  
USA  
E-MAIL: hsing@stat.ohio-state.edu

DEPARTMENT OF MATHEMATICAL STATISTICS  
CHALMERS UNIVERSITY OF TECHNOLOGY  
SE-41296 GOTHENBERG  
SWEDEN  
E-MAIL: rootzen@math.chalmers.se